\documentclass[12pt,leqno]{article}
\usepackage{amsmath,amssymb}
\begin{document}
\newtheorem{theorem}{Theorem}[section]
%\label %\ref
\newtheorem{prop}{Proposition}[section]
\newtheorem{defin}{Definition}[section]
\newtheorem{rem}{Remark}[section]
\newtheorem{example}{Example}[section]
\newtheorem{corol}{Corollary}[section]
\newtheorem{lemma}{Lemma}[section]
%\label{} %\ref{}
\title{Coupling Poisson and Jacobi structures on foliated manifolds}
\author{{\small by}\vspace{2mm}\\Izu Vaisman}
\date{}
\maketitle
{\def\thefootnote{*}\footnotetext[1]%
{{\it 2000 Mathematics Subject Classification: 53D17}.
Coupling Poisson structure. Coupling Jacobi structure.}}
\begin{center} \begin{minipage}{12cm}
A{\footnotesize BSTRACT. Let $M$ be a differentiable manifold
endowed with a foliation $\mathcal{F}$. A Poisson structure $P$ on
$M$ is $\mathcal{F}$-coupling if $\sharp_P(ann(T\mathcal{F}))$ is
a normal bundle of the foliation. This notion extends Sternberg's
coupling symplectic form of a particle in a Yang-Mills field
\cite{St}. In the present paper we extend Vorobiev's theory of
coupling Poisson structures \cite{Vor} from fiber bundles to
foliated manifolds and give simpler proofs of Vorobiev's existence
and equivalence theorems of coupling Poisson structures on duals
of kernels of transitive Lie algebroids over symplectic manifolds.
Then we discuss the extension of the coupling condition to Jacobi
structures on foliated manifolds.}
\end{minipage}
\end{center}
\vspace{5mm} \noindent The symplectic structure that describes the
coupling of a particle and a field was discovered by S. Sternberg
\cite{St} (for more details, see \cite{GLS}). In \cite{Vor}
Sternberg's symplectic form is extended to a general notion of a
coupling Poisson structure on a fiber bundle and the coupling
Poisson structures are used to get information about any Poisson
structure in the neighborhood of a symplectic leaf. In the present
paper we reprove Vorobiev's results in the context of foliated
manifolds and we extend the coupling condition to Jacobi
structures. In particular, we give simpler proofs for Vorobiev's
existence and equivalence theorems of coupling Poisson structures
on duals of kernels of transitive Lie algebroids over symplectic
manifolds. The motivation for studying Poisson and Jacobi
structures on a foliated manifold comes from the fact that a
foliated manifold may play the role of the phase space of a
physical system with gauge parameters, the latter being the
coordinates along the leaves.

{\it Acknowledgement}. The author is grateful to Yurii Vorobiev
for useful remarks during the work on this paper.
%begin{center} %\section %\end{center}
\section{Preliminaries}  Our framework is the $C^\infty$ category.
Let $M^m$ be an $m$-dimensional differentiable manifold and
$\mathcal{F}^p$ a regular foliation with $p$-dimensional leaves on
$M$ ($0\leq p\leq m$). We denote by $\Gamma$ spaces of global
cross sections of vector bundles and by
\begin{equation} \label{multivect} \mathcal{V}^k(M)=
\Gamma\wedge^kTM, \;\;\Omega^k(M)= \Gamma\wedge^kT^*M \hspace{5mm}
(k=1,\ldots,m) \end{equation} the spaces of multivector fields and
differential forms, respectively.

Furthermore, let us choose a $q$-dimensional, complementary
subbundle $H$ of $F=T\mathcal{F}$ $(q+p=m)$ i.e.,
\begin{equation} \label{descompTM} TM=H\oplus F.\end{equation}
$H$ is called a {\it normal bundle} of the foliation and we will
denote by $\pi_H,\pi_F$ the projections of $TM$ onto the
subbundles $H,F$, respectively.

The decomposition (\ref{descompTM}) induces a natural bigrading of
multivector fields and differential forms, and a corresponding
decomposition of the exterior differential
\begin{equation}\label{descompd}
d=d'_{1,0}+d''_{0,1}+\partial_{2,-1}, \end{equation} where the
indices denote the bidegree of the components in the direction of
$H,F$, respectively  \cite{V4}.

In this section, we give some computation formulas for
Schouten-Nijenhuis brackets (e.g., see \cite{V1}) on an arbitrary
differentiable manifold and on a foliated manifold with a fixed
decomposition (\ref{descompTM}).
\begin{prop}\label{lemaPP} For any bivector field
$P\in\mathcal{V}^2(M)$ one has
\begin{equation}\label{noulPP} [P,P](\alpha,\beta,\gamma)=
2[d\gamma(\sharp_P\alpha,\sharp_P\beta)-(L_{\sharp_P\gamma}P)(\alpha,\beta)],
\;\;\alpha,\beta,\gamma\in\Omega^1(M),
\end{equation}  where $L$ denotes the Lie derivative and
$\sharp_P:T^*M\rightarrow TM$ is defined by
$\beta(\sharp_P\alpha)=P(\alpha,\beta)$.
\end{prop} \noindent{\bf
Proof.} It would suffice to notice that the two sides of
(\ref{noulPP}) are tensor fields, then check the formula for
$\alpha,\beta,\gamma$ equal to differentials of local coordinates
on $M$. Instead, in order to recall some more very important
formulas, we will proceed as follows. For the {\it bracket of
$1$-forms}
\begin{equation}\label{bracket1} \{\alpha,\beta\}_P=
i(\sharp_P\alpha)d\beta-i(\sharp_P\beta)d\alpha +
d(P(\alpha,\beta)), \end{equation} one has the Gelfand-Dorfman
formula \cite{GD}
\begin{equation}\label{formula1GD} [P,P](\alpha,\beta,\gamma)=
2\{\gamma(\sharp_P\{\alpha,\beta\}_P-[\sharp_P\alpha,\sharp_P\beta])\}.
\end{equation} If expression (\ref{bracket1}) is inserted in
(\ref{formula1GD}), a computation that uses the classical relation
between the operators $d,i,L$ leads to (\ref{noulPP}). Q.e.d.

By polarizing formula (\ref{noulPP}) we get \begin{corol}
\label{PPpolar} For any two bivector fields
$P_1,P_2\in\mathcal{V}^2(M)$ one has
\begin{equation}\label{polarizare} [P_1,P_2](\alpha,\beta,\gamma)
=d\gamma(\sharp_{P_1}\alpha,\sharp_{P_2}\beta) +
d\gamma(\sharp_{P_2}\alpha,\sharp_{P_1}\beta)\end{equation} $$-
(L_{\sharp_{P_1}\gamma}P_2)(\alpha,\beta) -
(L_{\sharp_{P_2}\gamma}P_1)(\alpha,\beta).$$
\end{corol}

\begin{corol} \label{noilecond} The bivector field
$P\in\mathcal{V}^2(M)$ satisfies the Poisson condition $[P,P]=0$
iff \begin{equation}\label{nouacondP}
(L_{\sharp_P\gamma}P)(\alpha,\beta)=
d\gamma(\sharp_P\alpha,\sharp_P\beta).
\end{equation}
\end{corol}
\begin{rem} \label{obsalgLie}
{\rm In \cite{GD}, one also finds the following formula:
\begin{equation} \label{formula2GD}
\sum_{Cycl(\alpha,\beta,\gamma)}
<\{\{\alpha,\beta\},\gamma\},X>=[P,L_XP](\alpha,\beta,\gamma)
\end{equation} $$+\frac{1}{2}\sum_{Cycl(\alpha,\beta,\gamma)}
[P,P](\alpha,\beta,d<\gamma,X>),\hspace{5mm}\alpha,\beta,\gamma\in
\Omega^1(M),\;X\in\mathcal{V}^1(M).$$ Formulas (\ref{formula1GD})
and (\ref{formula2GD}) show that the cotangent bundle $T^*M$ of a
Poisson manifold $(M,P)$ with the bracket (\ref{bracket1}) is a
Lie algebroid of anchor $\sharp_P$, a fact which is of fundamental
importance in Poisson geometry.} \end{rem}

Now, assume the bivector field $P\in\mathcal{V}^2(M)$ is {\it
regular}, i.e., $s=rank\,P=const.$, and chose a decomposition
\begin{equation}\label{descompPregulat}
TM=E\oplus D,\;T^*M=E^*\oplus D^*\;\;(E^*=ann\,D,\;D^*=ann\,E),
\end{equation} where $D=im\,\sharp_P$ and we have corresponding natural
projections $p_E,p_D$. Then, we get an isomorphism
$\sharp_P:ann\,E\rightarrow D$ (by $ann$ we denote the annihilator
of a vector bundle), with an inverse $-\flat_P:D\rightarrow
ann\,E$, and there exists a well defined differential $2$-form of
rank $s$, $\theta\in\Gamma\wedge^2D^*$, defined by
\begin{equation} \label{2formaptS}
\theta(X,Y)=P[\flat_P(p_DX),\flat_P(p_DY)].\end{equation}
Conversely, (\ref{2formaptS}) allows us to reconstruct $P$ from
$\theta$, such that $ker\,\sharp_P=ann\,D$. We will say that
$\theta$ is {\it equivalent to $P$ modulo $E$}.
\begin{prop}\label{SScutheta} For a regular bivector field $P$,
$[P,P](\alpha,\beta,\gamma)=0$ if at least two of the arguments
belongs to $ann\,D$, and
\begin{equation} \label{SSregulat} \begin{array}{l}
[P,P](\alpha,\beta,\gamma)= 2\gamma([\sharp_P\alpha,\sharp_P\beta,
\hspace{5mm}\alpha,\beta\in ann\,E,\,\gamma\in ann\,D
\vspace{2mm}\\

[P,P](\alpha,\beta,\gamma)=
2d\theta(\sharp_P\alpha,\sharp_P\beta,\sharp_P\gamma),
\hspace{5mm}\alpha,\beta,\gamma\in ann\,E.
\end{array}\end{equation}
\end{prop} \noindent {\bf Proof.} We use the decomposition
(\ref{descompPregulat}). If not all the arguments are in $ann\,E$
the result immediately follows from either (\ref{noulPP}) or the
following Gelfand-Dorfman expression of the Schouten-Nijenhuis
bracket of two bivector fields \cite{GD}
\begin{equation}\label{NJGD}
[P,P](\alpha,\beta,\gamma)=2\sum_{Cycl(\alpha,\beta,\gamma)}
<\gamma,\sharp_P(L_{\sharp_P\alpha}\beta)>.\end{equation}  For
arguments in $ann\,E$, we have
$\alpha=\flat_PX,\beta=\flat_PY,\gamma=\flat_PZ$, with
$X,Y,Z\in\Gamma D$, and (\ref{NJGD}) yields
$$[P,P](\alpha,\beta,\gamma)=-2\sum_{Cycl(\alpha,\beta,\gamma)}
<L_X\beta,Z>$$ $$
=2\sum_{Cycl(X,Y,Z)}\{X(\theta(Y,Z))-\theta([X,Y],Z)\}.$$ Q.e.d.
\begin{corol}\label{strPreg}The regular bivector field $P$ is
Poisson iff the distribution $D=im\,\sharp_P$ is involutive and
there exists a decomposition $TM=E\oplus D$ such that the
equivalent to $P$, mod. $E$, $2$-form $\theta$ satisfies the
condition $d\theta|_{D}=0$, and, then, the same holds for any such
decomposition of $TM$.
\end{corol}

In the case of a foliated manifold $(M,\mathcal{F})$ with a normal
bundle $H$, if $P\in\mathcal{V}^2(M)$ is an arbitrary bivector
field, it has a decomposition
\begin{equation}\label{descompP}
P=P'_{2,0}+\bar P_{1,1}+P''_{0,2}, \end{equation} and we can
compute the corresponding decomposition of the Schouten-Nijenhuis
bracket $[P,P]$ by applying formula (\ref{noulPP}) to $1$-forms
$\alpha,\beta,\gamma\in\Omega^{1,0}(M)$ and
$\lambda,\mu,\nu\in\Omega^{0,1}(M)$. The results are contained in
the following formulas
\begin{equation}\label{descompSS1}\begin{array}{l}
[P',P']_{3,0}(\alpha,\beta,\gamma)=
2[d'\gamma(\sharp_{P'}\alpha,\sharp_{P'}\beta)
-(L_{\sharp_{P'}\gamma}P')(\alpha,\beta)],\vspace{2mm}\\

[P',P']_{2,1}(\alpha,\beta,\lambda) =
2\partial\lambda(\sharp_{P'}\alpha,\sharp_{P'}\beta) =
-2\lambda([\sharp_{P'}\alpha,\sharp_{P'}\beta]),
\vspace{2mm}\\

[P',P']_{1,2}(\alpha,\lambda,\mu)=0,\;\;
[P',P']_{0,3}(\lambda,\mu,\nu)=0,\end{array} \end{equation}
\begin{equation}\label{descompSS2} \begin{array}{l}
[\bar P ,\bar P]_{3,0}(\alpha,\beta,\gamma)=0,\vspace{2mm}\\

[\bar P ,\bar P]_{2,1}(\alpha,\beta,\lambda)
=2[d''\lambda(\sharp_{\bar P}\alpha,\sharp_{\bar P}\beta)
-(L_{\sharp_{\bar P}\lambda}\bar P)(\alpha,\beta)],\vspace{2mm}\\

[\bar P ,\bar P]_{1,2}(\alpha,\lambda,\mu) =-2[d'\mu(\sharp_{\bar
P}\lambda,\sharp_{\bar P}\alpha) +(L_{\sharp_{\bar P}\mu}\bar
P)(\alpha,\lambda)],\vspace{2mm}\\

[\bar P ,\bar P]_{0,3}(\lambda,\mu,\nu)
=2[\partial\nu(\sharp_{\bar P}\lambda,\sharp_{\bar P}\mu)
-(L_{\sharp_{\bar P}\nu}\bar P)(\lambda,\mu)],\vspace{2mm}
\end{array}\end{equation}
\begin{equation}\label{descompSS3}
\begin{array}{l} [P'',P'']_{3,0}(\alpha,\beta,\gamma)=0,\;\;
[P'',P'']_{2,1}(\alpha,\beta,\lambda)=0,\vspace{2mm}\\

[P'',P'']_{1,2}(\alpha,\lambda,\mu)=0,\vspace{2mm}\\

[P'',P'']_{0,3}(\lambda,\mu,\nu)
=2[d''\nu(\sharp_{P''}\lambda,\sharp_{P''}\mu)
-(L_{\sharp_{P''}\nu}P'')(\lambda,\mu)], \end{array}\end{equation}
\begin{equation}\label{descompSS4}
\begin{array}{l}\hspace*{-10mm} [P',\bar P]_{3,0}(\alpha,\beta,\gamma)
=d''\gamma(\sharp_{P'}\alpha,\sharp_{\bar P}\beta) -
d''\gamma(\sharp_{P'}\beta,\sharp_{\bar P}\alpha)\vspace{2mm}\\

\hspace*{-10mm}\hspace*{3cm} -(L_{\sharp_{P'}\gamma}\bar
P)(\alpha,\beta) - (L_{\sharp_{\bar
P}\gamma}P')(\alpha,\beta),\vspace{2mm}\\

\hspace*{-10mm}[P',\bar P]_{2,1}(\alpha,\beta,\lambda)
=d'\lambda(\sharp_{P'}\alpha,\sharp_{\bar P}\beta) -
d'\lambda(\sharp_{P'}\beta,\sharp_{\bar P}\alpha) \vspace{2mm}\\

\hspace*{-10mm}\hspace*{3cm}-(L_{\sharp_{\bar
P}\lambda}P')(\alpha,\beta),\vspace{2mm}\\

\hspace*{-10mm}[P',\bar P]_{1,2}(\alpha,\lambda,\mu)
=\partial\mu(\sharp_{P'}\alpha,\sharp_{\bar P}\lambda) -
(L_{\sharp_{\bar P}\mu}P')(\alpha,\lambda)\vspace{2mm}\\

\hspace*{-1cm}\hspace*{1cm}=-(L_{\sharp_{P'}\alpha}\bar P)
(\lambda,\mu), \;\; [P',\bar
P]_{0,3}(\lambda,\mu,\nu)=0,\end{array}
\end{equation}
\begin{equation}\label{descompSS5}
\begin{array}{l}\hspace*{-4mm}
[P',P'']_{3,0}(\alpha,\beta,\gamma)=0,\;\;
[P',P'']_{0,3}(\lambda,\mu,\nu)=0,
\vspace{2mm}\\

\hspace*{-4mm}[P',P'']_{2,1}(\alpha,\beta,\lambda)=
-(L_{\sharp_{P''}\lambda}P')(\alpha,\beta)
\vspace{2mm}\\

\hspace*{-4mm}[P',P'']_{1,2}(\alpha,\lambda,\mu)=
d'\mu(\sharp_{P'}\alpha,\sharp_{P''}\lambda)
-(L_{\sharp_{P''}\mu}P')(\alpha,\lambda)
\vspace{2mm}\\

\hspace*{-4mm}\hspace*{3cm}
=(L_{\sharp_{P'}\alpha}P'')(\lambda,\mu),
\end{array}
\end{equation}
\begin{equation}\label{descompSS6}
\begin{array}{l}\hspace*{-5mm}
[\bar P,P'']_{3,0}(\alpha,\beta,\gamma)=0,\;\;[\bar
P,P'']_{2,1}(\alpha,\beta,\lambda)=0, \vspace{2mm}\\

\hspace*{-5mm}[\bar P,P'']_{1,2}(\alpha,\lambda,\mu)=
d''\mu(\sharp_{\bar P}\alpha,\sharp_{P''}\lambda)
-(L_{\sharp_{\bar P}\mu}P'')(\alpha,\lambda) \vspace{2mm}\\

\hspace*{-5mm}\hspace*{3cm} -(L_{\sharp_{P''}\mu}\bar
P)(\alpha,\lambda),\vspace{2mm}\\

\hspace*{-5mm}[\bar P,P'']_{0,3}(\lambda,\mu,\nu)
=d'\nu(\sharp_{\bar P}\lambda, \sharp_{P''}\mu) - d'\nu(
\sharp_{\bar
P}\mu,\sharp_{P''}\lambda)\vspace{2mm}\\

\hspace*{-5mm}\hspace*{3cm} -(L_{\sharp_{\bar
P}\nu}P'')(\lambda,\mu)-(L_{\sharp_{P''}\nu}\bar P)(\lambda,\mu).
\end{array} \end{equation}

Finally, we indicate a way to compute the bigraded components of
any Schouten-Nijenhuis bracket. Let us recall the following
general formula of Lichnerowicz (e.g., see \cite{V1}):
\begin{equation} \label{Lichne} i([P,Q])\varphi= (-1)^q(p+1)i(P)
d(i(Q)\varphi)\end{equation} $$+(-1)^pi(Q)d(i(P)\varphi)
-i(P\wedge Q)d\varphi,$$ where $P\in\mathcal{V}^p(M),
Q\in\mathcal{V}^q(M)$, and $\varphi\in\Omega^{p+q-1}(M)$. The
operator $i$ is the {\it interior product}. If $g$ is an arbitrary
Riemannian metric on $M$, $i(P)$ is the transposed operator of the
wedge product by $\flat_gP$ \cite{Lch}, whence, it follows that
\begin{equation} \label{idewedge}i(P\wedge Q)=i(Q)\circ
i(P).\end{equation}

Now, if $P$ is of bidegree $(a,b)$ $(a+b=p)$ and $Q$ is of
bidegree $(h,k)$ $(h+k=q)$, the component $[P,Q]^{u,v}$
$(u+v=p+q-1)$ is provided by using (\ref{Lichne}) for
$\varphi\in\Omega^{u,v}(M)$. With (\ref{descompd}), we see that
the only possibilities to get non-zero components correspond to
the replacement of $d$ by $d',d'',\partial$ in (\ref{Lichne}),
which leads to the cases
\begin{equation} \label{compinS-Nij} u= a+h-1,
v=b+k;\;u=a+h,v=b+k-1; \end{equation}
$$u=a+h-2,\;v=b+k+1.$$ All the other components vanish because of
degree incompatibility.
\section{Leaf-tangent Poisson structures} In this section we
discuss Poisson structures along the leaves of a foliation since
these will be a basic ingredient of the coupling Poisson
structures. The notation is that of Section 1.
\begin{defin}
\label{strpefoi} {\rm The Poisson structure defined by the bivector
field $P\in\mathcal{V}^2(M)$ is {\it leaf-tangent} to
$\mathcal{F}$ if its symplectic leaves are submanifolds of the
leaves of $\mathcal{F}$, equivalently, if the leaves of
$\mathcal{F}$ are Poisson submanifolds of $(M,P)$
\cite{V1}.}
\end{defin}

Obviously, $P$ is $\mathcal{F}$-leaf-tangent iff
$P\in\Gamma\wedge^2F$ $(F=T\mathcal{F})$, which is equivalent with
the fact that, for all the choices of a normal bundle $H$
satisfying (\ref{descompTM}), the decomposition (\ref{descompP})
of $P$ satisfies the conditions
\begin{equation} \label{condPpefoi} P'_{1,1}=0,\;\;\bar P_{1,2}=0.
\end{equation}
\begin{prop} \label{tangentfoi} A bivector field
$P\in\Gamma\wedge^2F$ is a Poisson bivector field on $M$ iff
\begin{equation} \label{SS3ptP}
d''\nu(\sharp_{P}\lambda,\sharp_{P}\mu)
-(L_{\sharp_{P}\nu}P)(\lambda,\mu)=0, \hspace{5mm}
\forall\lambda,\mu,\nu\in\Omega^{0,1}(M),
\end{equation} equivalently,
iff the restrictions of $P$ to the leaves are Poisson bivector
fields of the leaves.
\end{prop} \noindent {\bf Proof.}
$P$ is Poisson iff $[P,P]=0$ and, by
(\ref{descompSS1})-(\ref{descompSS6}), this condition reduces to
(\ref{SS3ptP}). Q.e.d.

Therefore, a leaf-tangent Poisson structure may be identified with
a family $P_L$ of Poisson structures of the leaves $L$ of
$\mathcal{F}$ such that, $\forall f,g\in C^\infty(M)$, the
function
$$\{f,g\}(x)=\{f|_{L(x)},g|_{L(x)}\}_{P|_{L(x)}}(x),$$ where $x\in
M$ and $L(x)$ is the leaf of $\mathcal{F}$ through $x$, also
belongs to $C^\infty(M)$.

It is worthwhile noticing that, if $P$ is
$\mathcal{F}$-leaf-tangent, $\mathcal{F}$ may be seen as a {\it
regularizing foliation} of the symplectic foliation $\mathcal{S}$
of $P$ (usually, $\mathcal{S}$ has singular points). This idea
leads to a global, numerical invariant of a Poisson structure on a
manifold $M^m$ , the {\it regularizing dimension}, to be defined
as the smallest possible dimension $p\leq m$ of a regularizing
foliation.
\begin{example}\label{ex3-d} {\rm Put
$\mathbb{R}^m=\mathbb{R}^3\times\mathbb{R}^{m-3}$ and take the
Poisson structure $P$ defined by the Lie-Poisson structure of the
factor $\mathbb{R}^3$ seen as the Lie algebra $so(3)$ (e.g., \cite
{V1}). $P$ is leaf-tangent to the foliation defined by the factor
$\mathbb{R}^3$ of $\mathbb{R}^m$, and the regularizing dimension
of $P$ is $3$.}\end{example}
\begin{example} \label{Heisenberg}
{\rm Consider the compact nilmanifold
$M(1,n)=\Gamma(1,n)\backslash H(1,n)$ where
\begin{equation} \label{eqHeisenberg}
H(1,n)=\left\{\left(\begin{array}{ccc} Id_n&X&Z\\ 0&1&y\\0&0&1
\end{array}\right)\,/\,X,Z\in \mathbb{R}^n,\,y\in
\mathbb{R}\right\}\end{equation} is the generalized Heisenberg
group, and $\Gamma(1,n)$ is the subgroup of matrices with integer
entries. $M(1,n)$ has a natural atlas with the transition
functions
\begin{equation} \label{trHeis} \tilde x^i=x^i+a^i,\,\tilde
y=y+b,\, \tilde z^i=z^i+a^iy+c^i,\end{equation} where $x^i,z^i$
$(i=1,...,n)$ are the entries of $X,Z$, respectively, and
$a^i,b,c^i$ are integers, and it is parallelizable by the global
vector fields
\begin{equation}\label{paralHeis}
\frac{\partial}{\partial x^i}, \frac{\partial}{\partial
y}+\sum_{i=1}^px^i\frac{\partial}{\partial
z^i},\frac{\partial}{\partial z^i}.\end{equation}

It follows that \begin{equation} \label{PpeHeis} P=\sin{(2\pi
y)}\frac{\partial}{\partial x^1}\wedge\frac{\partial}{\partial
z^1}\end{equation} yields a Poisson structure of $M(1,n)$ which is
leaf-tangent to the $(n+1)$-dimensional foliation
\begin{equation}\label{fol2exHeis} y=const.,\; x^t=const.
\hspace{5mm}(t=2,...,n). \end{equation}

But, $P$ is also leaf-tangent to the $2$-dimensional foliation
$span\{\partial/\partial x^1,
\partial/\partial z^1\}$, and the
regularizing dimension of $P$ is $2$.}
\end{example}
\begin{example} \label{fibratePoisson} {\rm \cite{Vor} Let
$p:\mathbb{G}^*\rightarrow B$ be a bundle of Lie coalgebras (i.e.,
the dual of a bundle $p:\mathbb{G}\rightarrow B$ of Lie algebras)
over a manifold $B$. Then the Lie-Poisson structures of the fibers
yield a leaf-tangent Poisson structure $\mathbb{L}$ of
$\mathbb{G}^*$. If $(x^i)$ are local coordinates on $B$ and
$(y_a)$ are linear coordinates along the fibers of $\mathbb{G}^*$,
one has \cite{V1}
\begin{equation}\label{PpemathcalG}
\mathbb{L}=\frac{1}{2}\alpha_{ab}^c(x)y_c \frac{\partial}{\partial
y_a}\wedge \frac{\partial}{\partial y_b}, \end{equation} where the
Einstein summation convention is used and $\alpha_{ab}^c(x)$ are
the structural constants of the corresponding fibers of
$\mathbb{G}$. The invariant expression of $\mathbb{L}$ on
differentials at $z\in\mathbb{G}^*$ of fiberwise linear functions
on $\mathbb{G}^*$, seen as elements $X,Y\in\mathbb{G}_{p(z)}$, is
\begin{equation} \label{LiePinvariant} \mathbb{L}_z(X,Y)=
<z,C_{p(z)}(X,Y)>,\end{equation} where the ``tensor field"
$C\in\Gamma[(\wedge^2\mathbb{G}^*)\otimes (\mathbb{G})]$ is
defined by
\begin{equation} \label{tensorC} C_{p(z)}(X,Y)=
[X,Y]_{\mathbb{G}_{p(z)}}.\end{equation}}
\end{example}

The properties of a leaf-tangent Poisson structure reflect
corresponding properties along the leaves. An interesting
situation appears for the Poisson cohomology.
\begin{prop} \label{descompsigma} Let $P$ be a leaf-tangent
Poisson structure on $(M,\mathcal{F})$ and let $\sigma$ be the
Lichnerowicz coboundary operator of $P$. Then, for every normal
bundle $H$ of $\mathcal{F}$, one has a decomposition
\begin{equation} \label{descompsig''}
\sigma=\sigma'_{-1,2}+\sigma''_{0,1} \end{equation}
where the terms are homogeneous components of $\sigma$  of
indicated bidegree and
\begin{equation}\label{bicomplex} \sigma'^2=0,\;
\sigma''^2=0,\;\sigma'\circ\sigma''+ \sigma''\circ\sigma'=0. \end{equation}
\end{prop} \noindent{\bf Proof.} The Lichnerowicz coboundary operator is
$\sigma Q=-[P,Q]$, $\forall Q\in\mathcal{V}^*(M)$ (e.g.,
\cite{V1}). From (\ref{compinS-Nij}), since $P$ is homogeneous of
bidegree $(0,2)$, it follows that $\sigma$ may have only
components of bidegree $(-1,2)$, $(0,1)$ and $(-2,3)$. We will
prove that the $(-2,3)$-component of $\sigma$ vanishes.

Indeed, $\sigma_{-2,3}$ is obtained by computing
$i([P,Q])\varphi$, $Q\in\mathcal{V}^{h,k}(M)$,
$\varphi\in\Omega^{h-2,k+3}$ via (\ref{Lichne}). Since in this
case $i(Q)\varphi=0$ (it should be of bidegree $(-2,3)$), using
(\ref{idewedge}) we get
\begin{equation} \label{calculpti}
i([P,Q])\varphi=i(Q)[di(P)\varphi- i(P)d\varphi],\end{equation}
where, in fact, only the component $\partial$ of $d$ may bring a
non zero contribution.

Now, let us evaluate a Lie derivative $L_Y\psi$ where $Y\in\Gamma
F$ and $\psi\in
\Omega^{s,t}(M)$ for arguments
$X_1,...,X_u\in\mathcal{V}_{pr}^{1,0}(M),\;
Y_1,...,Y_v\in\mathcal{V}^{0,1}(M)$, where
$\mathcal{V}^{1,0}_{pr}(M)$ is the subspace of the $(1,0)$-vector
fields that project onto the space of leaves of $\mathcal{F}$,
which means that, $\forall Y\in\Gamma F$, $[Y,X_i]\in \Gamma F$
$(i=1,...,u)$. (It is enough to use only this kind of arguments
$X$ because $L_Y\psi$ is a tensor, hence, its values at each point
depend on the values of the arguments at that point only.) We have
$$(L_Y\psi)(X_1,...,X_u,Y_1,...,Y_v)=Y(\psi(X_1,...,X_u,Y_1,...,Y_v))$$
$$-\sum_{i=1}^u\psi(X_1,...,X_{i-1},[Y,X_i],X_{i+1},...,X_u,Y_1,...,Y_v)$$
$$-\sum_{j=1}^v\psi(X_1,...,X_u,Y_1,...,Y_{j-1},[Y,Y_j],Y_{j+1},...,Y_v)=0,$$
and it follows that $L_Y\psi$ may have only components of bidegree
equal either to $(s,t)$ or to $(s-1,t+1)$.

From the previous remark, it follows that, $\forall
Q\in\mathcal{V}^{h,k}(M)$, $\forall\psi\in
\Omega^{h-2,k+3}\oplus\Omega^{h-2,k+2}$ and $\forall Y\in\Gamma F$, one has
\begin{equation} \label{eqauxiliar} i(Q)L_Y\psi=
i(Q)[di(Y)+i(Y)d]\psi=0. \end{equation}

Furthermore, for $\varphi\in\Omega^{h-2,k+3}$, (\ref{idewedge})
and (\ref{eqauxiliar}) imply
$$i(Q)[i(Y_1\wedge Y_2)d-di(Y_1\wedge Y_2)]\varphi
=i(Y_2\wedge Q)i(Y_1)d\varphi-i(Q)di(Y_2)i(Y_1)\varphi$$
$$=-i(Y_2\wedge Q)di(Y_1)\varphi+i(Q)i(Y_2)di(Y_1)\varphi=0.$$
Since $P$ is spanned over $\mathbb{R}$ by wedge products
$Y_1\wedge Y_2$ of vector fields tangent to $\mathcal{F}$, the
expression given by (\ref{calculpti}) vanishes.

The properties (\ref{bicomplex}) follow from $\sigma^2=0$. Q.e.d.

Even though we have (\ref{bicomplex}), the degrees are not right
for a double cochain complex. However,
$(\mathcal{W}^{h,k}=\mathcal{V}^{k,h}(M),\sigma)$ is a double,
semipositive, cochain complex, and the cohomology of such a
complex is the limit of a spectral sequence
\cite{V4}. More exactly, we have
\begin{prop} \label{propspectral} Let $P$ be a
leaf-tangent Poisson structure on $(M,\mathcal{F})$ and denote by
$H_{LP}^h(\mathcal{F},P)$ the cohomology spaces of the cochain
complex $(\mathcal{V}^{0,*},\sigma'')$. Then, the Poisson
cohomology of $P$ is the limit of a spectral sequence
$(E^{hk}_r,d_r)$ where
\begin{equation}
\label{Edoi} E_2^{hk}=
\mathcal{V}^{k,0}(M)\otimes_\mathbb{R} H_{LP}^h(\mathcal{F},P)\end{equation}
and $d_2$ is induced by the operator $\sigma'$. \end{prop}
\noindent {\bf Proof.} The spaces
$$\mathcal{W}_l(M)=\oplus_{k\geq l}\oplus_h\mathcal{V}^{hk}(M)$$
yield a regular filtration of the Lichnerowicz-Poisson complex
$(\mathcal{V}(M),\sigma)$ and the required spectral sequence is
the spectral sequence $(E^{hk}_r,d_r)$ defined by this filtration.
From the definition of a spectral sequence we get
$$E_0^{kh}=\mathcal{V}^{h,k}(M),\;\;d_0=0.$$
The cohomology of $E_0$ yields
$$E_1^{kh}=\mathcal{V}^{h,k}(M),\;\;d_1=\tau''.$$
Finally, the cohomology spaces of the complex $E_1$ are the spaces
given by formula (\ref{Edoi}) and $d_2$  is induced by $\tau'$.
Q.e.d.

\begin{rem} \label{obscohF} {\rm The leaf-tangent Poisson bivector field
$P\in\Gamma\wedge^2F$ may be seen as a Poisson bivector of the Lie
algebroid $F$ defined by the foliation. As such, it induces a Lie
algebroid structure on the dual bundle $F^*$, and
$H_{LP}^h(\mathcal{F},P)$ are the cohomology spaces of this Lie
algebroid (e.g.,
\cite{MX}).} \end{rem}
\section{Coupling Poisson structures}
In this section we present the general results concerning coupling
Poisson structures on foliated manifolds. Again, the notation is
that of Section 1, we consider the manifold $M$, the foliation
$\mathcal{F}$, the normal bundle $H$, and the bivector field $P$
written under the form (\ref{descompP}).
\begin{defin} \label{compatalg}
{\rm The bivector field $P$ is $\mathcal{F}$-{\it almost coupling
via} $H$ if, $\forall x\in M$,
\begin{equation} \label{compatalgVorob} \sharp_P(ann\,F_x)\subseteq
H_x. \end{equation} ii) The bivector field $P$ is
$\mathcal{F}$-{\it coupling} if $\sharp_P(ann\,F)$ is a normal
bundle $H$ of $\mathcal{F}$. In both cases, if $P$ is Poisson, the
coupling property is attributed to the Poisson structure.}
\end{defin}

With (\ref{descompP}), it follows that the almost coupling
condition is equivalent with the vanishing of the
$(1,1)$-component $\bar P=0$, hence, also with
$\sharp_P(ann\,H)\subseteq F$. The coupling condition is
equivalent with \begin{equation} \label{rakcondition}
dim(\sharp_P(ann\,F))=q,
\end{equation} the codimension of $\mathcal{F}$. The first term of
(\ref{rakcondition}) is the rank of the term $P'$ of
(\ref{descompP}) for any choice of $H$, hence, coupling may exist
only if $q$ is even. In the case of a coupling field $P$, we
always shall take $H=\sharp_P(ann\,F)$.
\begin{example} \label{ex1coupling} {\rm Consider again the
manifold $M(1,n)$ of Example \ref{Heisenberg}. With the notation
of that example, define the foliation $\mathcal{F}$ by the tangent
distribution \begin{equation} \label{eqex1coupling} F=span\left\{
\frac{\partial}{\partial x^1}, \frac{\partial}{\partial z^1},
\frac{\partial}{\partial y} + \sum_{i=1}^n
x^i\frac{\partial}{\partial z^i}\right\}
\end{equation} and the Poisson bivector field
\begin{equation} \label{eq2ex1coupling}
P=\sum_{t=2}^n \frac{\partial}{\partial x^t}\wedge
\frac{\partial}{\partial z^t} + \sin{2\pi y}
\frac{\partial}{\partial x^1}\wedge \frac{\partial}{\partial z^1}.
\end{equation} Obviously, $P$ is
$\mathcal{F}$-coupling.} \end{example}
\begin{example} {\rm Let $P$ be an arbitrary Poisson bivector
field on the foliated manifold $(M,\mathcal{F})$. Assume that
there exists a symplectic leaf $S$ of $P$ which is embedded in $M$
and transversal to $\mathcal{F}$. (Generally, such a situation may
appear on an open subset $U$ of $M$.) Then condition
(\ref{rakcondition}) holds on $S$, therefore, also on some open
neighborhood $V$ of $S$, and $P$ is $\mathcal{F}$-coupling on $V$.
In particular, for any Poisson structure $P$ and any symplectic
leaf $S$ of $P$ which is embedded in $M$, there exists a tubular
neighborhood $V$ of $S$ where $P$ is coupling with respect to the
fibers of the tubular structure of $V$. This remark is due to
Vorobiev \cite{Vor}.}
\end{example}
\begin{prop} \label{compalgPois} An almost coupling bivector field
$P$ is Poisson iff
\begin{equation} \label{farabarP} \begin{array}{l}
 d'\gamma(\sharp_{P'}\alpha,\sharp_{P'}\beta)
-(L_{\sharp_{P'}\gamma}P')(\alpha,\beta)=0,\vspace{2mm}\\

(L_{\sharp_{P''}\lambda}P')(\alpha,\beta)+\lambda([\sharp_{P'}\alpha,
\sharp_{P'}\beta])=0,\vspace{2mm}\\

(L_{\sharp_{P'}\alpha}P'')(\lambda,\mu)=0,\vspace{2mm}\\

 d''\nu(\sharp_{P''}\lambda,\sharp_{P''}\mu)
-(L_{\sharp_{P''}\nu}P'')(\lambda,\mu)=0,
\end{array} \end{equation}
$\forall\alpha,\beta,\gamma\in\Omega^{1,0}(M)$,
$\forall\lambda,\mu,\nu\in\Omega^{0,1}(M)$.
\end{prop} \noindent {\bf Proof.}
Use (\ref{descompSS1})-(\ref{descompSS6}) to express $[P,P]=0$ for
\begin{equation} \label{Pcoupling} P=P'+P''.\end{equation} Q.e.d.

The last condition (\ref{farabarP}) is (\ref{SS3ptP}), again, and
it means that the component $P''$ is an $\mathcal{F}$-leaf-tangent
Poisson bivector field.

In the coupling case, the Poisson condition may be put in the form
given by Vorobiev \cite{Vor} for fiber bundles.
\begin{prop}\label{thVorobiev} A coupling Poisson structure of a
foliated manifold $(M,\mathcal{F})$ is equivalent with a triple
$(P'',H,\sigma)$, where $P''$ is a leaf-tangent Poisson structure,
$H$ is a normal bundle of $\mathcal{F}$ and $\sigma$ is a
non-degenerate cross section of $\wedge^2(ann\,F)$ such that
\begin{equation}\label{d'sigmazero}
d'\sigma=0,
\end{equation}
\begin{equation}\label{condcurb}
\sharp_{P''}\{d[\sigma(X,Y)]\}=-p_F[X,Y],\hspace{5mm}\forall
X,Y\in \mathcal{V}^{1,0}_{pr}(M),
\end{equation}
\begin{equation}\label{condLie} L_XP''=0,
\hspace{1cm}\forall X\in \mathcal{V}^{1,0}_{pr}(M).\end{equation}
\end{prop} \noindent {\bf Proof.} Let $P$ be a coupling Poisson
bivector field written under the form (\ref{Pcoupling}). The
coupling condition defines $H$ and $P''$, and the last condition
(\ref{farabarP}) is equivalent with the fact that $P''$ is a
Poisson bivector field. Furthermore, $P'$ is a regular bivector
field with $im\,\sharp_{P'}=H$, and we may write
$P'=\sharp_{P'}(\sigma)$ $=\sharp_P(\sigma)$, where
$\sigma\in\wedge^2(ann\,F)$ is equivalent with $P'$ mod. $F$ and
has the maximal rank. Accordingly, for the $1$-forms of bidegree
$(1,0)$ of (\ref{farabarP}) we may write
\begin{equation} \label{reprezcuXYZ} \alpha=\flat_\sigma X,\;
\beta=\flat_\sigma Y,\; \gamma=\flat_\sigma Z,\hspace{5mm}
X,Y,Z\in\mathcal{V}^{1,0}(M),\end{equation} where $X,Y,Z$ are
uniquely defined. With this representation, formulas
(\ref{noulPP}), (\ref{SSregulat}) show that the first condition
(\ref{farabarP}) becomes (\ref{d'sigmazero}).

Conditions (\ref{farabarP}) are tensorial. Thus, in particular,
the second condition (\ref{farabarP}) holds iff it holds for
$X,Y\in\mathcal{V}^{1,0}_{pr}(M)$ in (\ref{reprezcuXYZ}). In this
case, using the definition of the Lie derivative, we see that the
second condition (\ref{farabarP}) becomes (\ref{condcurb}).

Finally, the third condition (\ref{farabarP}) is equivalent to
(\ref{condLie}) because, for a projectable vector field $X$,
$L_XP''$ vanishes if at least one of its arguments is of bidegree
$(1,0)$.

For the converse, first notice that we may define a triple
$(P'',H,\sigma)$ for any coupling bivector field $P$. The $2$-form
$\sigma$ will be called a {\it coupling form} of the leaf-tangent
bivector field $P''$ and one has an isomorphism
$\flat_\sigma:H\rightarrow ann\,F$. Then, we can reconstruct $P$
from such a triple by the formula
\begin{equation}\label{Pdintriplu} P(\xi,\eta)=
\sigma(\flat_\sigma^{-1}\xi',\flat_\sigma^{-1}\eta')
+P''(\xi'',\eta''),\end{equation} where
$$\xi=\xi'+\xi'',\;\eta=\eta'+\eta'',\hspace{5mm}
\xi',\eta'\in ann\,F,\,\xi'',\eta''\in ann\,H.$$
If we reconstruct $P$ from $(P'',H,\sigma)$ that satisfy
(\ref{d'sigmazero}), (\ref{condcurb}), (\ref{condLie}), and where
$P''$ is Poisson, conditions (\ref{farabarP}) will hold. Q.e.d.
\begin{rem}\label{obssuplementara2} {\rm If the triple
$(P'',H,\sigma)$ satisfies all the conditions of Proposition
\ref{thVorobiev}, these conditions are also satisfied by any
triple $(P'',H,\sigma+\epsilon\tau)$ where
$\tau\in\wedge^2(ann\,F)$ is closed and $\epsilon\in\mathbb{R}$.
If $\epsilon$ is small enough, $\sigma+\epsilon\tau$ is non
degenerate on $H$, and the new triple also provides a coupling
Poisson structure.}\end{rem}
\begin{rem}\label{obssuplementara} {\rm
If the coupling Poisson tensor field $P$ of (\ref{Pcoupling}) is
defined by a symplectic form $\omega$ of a manifold $M$, $P''$ is
equivalent with a closed $2$-form $\theta$ of bidegree $(0,2)$
that defines a symplectic structure on each leaf of $\mathcal{F}$,
and the form $\sigma$ of the corresponding triple is a
corresponding coupling form in the sense of
\cite{GLS}.} \end{rem}

We also notice the following, rather clear, assertion
\begin{prop}\label{couplingequivalence} Two coupling bivector
fields $P_1,P_2$ are equivalent by a foliated diffeomorphism
$\Phi$ iff the corresponding triples are equivalent by $\Phi$.
\end{prop} \noindent {\bf Proof.} By equivalence we mean
that $\Phi_*(P_1)=P_2$, i.e.,
\begin{equation}\label{eqechivalP12}
\sharp_{P_2(\phi(x))}=\Phi_*(x)\circ\sharp_{P_1(x)}\circ\Phi^*(\Phi(x)),
\hspace{5mm}(x\in M),\end{equation}  and a diffeomorphism
$\Phi:M\rightarrow M$ is {\it foliated} if it preserves the
foliation $\mathcal{F}$, whence, $\Phi_*(F_x)=F_{\Phi(x)}$ and
$\Phi^*(ann\,F_{\Phi(x)})=ann\,F_x$. These properties, and
(\ref{eqechivalP12}), imply that, if $P_1$ is coupling, $P_2$
necessarily is coupling as well and it defines the normal bundle
$H_2=\Phi_*(H_1)$. Now, using the decompositions
$$TM=H_1\oplus F=H_2\oplus F$$ it follows that
$P''_2=\Phi_*(P''_1)$ and $\sigma_2=\Phi^*\sigma_1$. Conversely,
the equivalence of the triples implies that of the bivector fields
in view of formula (\ref{Pdintriplu}). Q.e.d.
\begin{rem}\label{echivalnefol} {\rm A Poisson equivalence of
coupling Poisson structures which is not a foliated diffeomorphism
does not preserve the elements of the associated triples.}
\end{rem}

The conditions indicated in Proposition \ref{thVorobiev} can be
expressed in a more detailed way if $P''$ is a regular bivector
field. \begin{prop}\label{Vorobregular} A regular coupling Poisson
structure of a foliated manifold $(M,\mathcal{F})$ is equivalent
with a class of objects $(S,G,H,\sigma,\theta)$, where $S,G,H$ are
vector subbundles of $TM$ such that: i) $S$ is integrable,
$F=G\oplus S$ and $TM=H\oplus F$; ii)
$\sigma\in\Gamma\wedge^2(ann\,F)$ is non degenerate and satisfies
the condition $d'\sigma=0$; iii)
$\theta\in\Gamma\wedge^2(ann(H\oplus G))$ is non degenerate and
satisfies the condition $d\theta|_S=0$; iv) for any projectable
vector field $X\in\Gamma H$ and any vector fields $U,V\in\Gamma S$
one has $[X,U]\in\Gamma S$ and $(L_X\theta)(U,V)=0$; v) for any
projectable vector fields $X,Y\in\Gamma H$, the projection
$p_F[X,Y]\in\Gamma S$ and $\flat_\theta(p_F[X,Y])=d[\sigma(X,Y)]$,
where $\flat_\theta$ is the isomorphism $S\rightarrow ann(H\oplus
G)$ defined by $\theta$.
\end{prop} \noindent {\bf Proof.}
Consider the associated triple $(P'',H,\sigma)$ that satisfies the
conditions of Proposition \ref{thVorobiev}, and ask $P''$ to have
a constant rank. Let $S$ be the tangent bundle of the symplectic
foliation of $P''$, let $G$ be an arbitrary tangent subbundle such
that $F=G\oplus S$, and let $\theta$ be the equivalent $2$-form of
$P''$, mod. $H\oplus G$. Then formulas (\ref{noulPP}),
(\ref{SSregulat}) and the integrability of $S$ show that the last
condition (\ref{farabarP}) is exactly $d\theta|_S=0$.

As indicated in the proof of Proposition \ref{thVorobiev},
(\ref{condLie}) holds iff it holds for  arguments
$\lambda,\mu\in\Gamma(ann\,H)$, which may be further decomposed as
$\lambda=\lambda_1+\lambda_2$, $\mu=\mu_1+\mu_2$, where
$\lambda_1,\mu_1\in\Gamma(ann(H\oplus S))$,
$\lambda_2,\mu_2\in\Gamma(ann(H\oplus G))$. Then,
$\sharp_{P''}\lambda_1=0, \sharp_{P''}\mu_1=0$, and
$\lambda_2=\flat_\theta U$, $\mu_2=\flat_\theta V$ with
$U,V\in\Gamma S$, and $\forall X\in \mathcal{V}^{1,0}_{pr}(M)$ we
get
\begin{equation}\label{eqajutatoare}
(L_X P'')(\lambda_1+\lambda_2,\mu_1+\mu_2)=
X(\theta(V,U))+\theta(U,[X,V])-\theta(V,[X,U])\end{equation}
$$+\lambda_1([X,V])
-\mu_1([X,U])=-(L_X\theta)(U,V)+\lambda_1([X,V])-\mu_1([X,U]).$$
Since $\lambda_1,\mu_1\in\Gamma(ann(H\oplus S))$ are arbitrary
arguments, the result of (\ref{eqajutatoare}) is zero iff
condition iv) of the proposition holds.

The equivalence between (\ref{condcurb}) and condition v) is
obvious.

The computations above also show that a system
$(S,G,H,\sigma,\theta)$ which satisfies the conditions of the
present proposition produces a triple $(P'',H,\sigma)$ that
satisfies the conditions of Proposition \ref{thVorobiev}. Two
systems of data $(S,G,H$, $\sigma,\theta)$ will be in the same
class, in the sense of the proposition, if they produce the same
bivector field $P''$. Q.e.d.

In particular, we get the characteristic conditions of \cite{GLS}
for the symplectic coupling forms:
\begin{corol} \label{symplecticcoupling} Let $(M,\mathcal{F})$
be a foliated manifold endowed with a symplectic form $\omega$.
Then $\omega$ is $\mathcal{F}$-coupling iff  the
$\omega$-orthogonal distribution $H$ of $F$ is a normal bundle of
$\mathcal{F}$ and one has a corresponding decomposition
\begin{equation} \label{descsymplcoup} \omega=\sigma_{2,0}+\theta_{0,2}
\end{equation} that satisfies the following conditions
i) the form $\theta$ is symplectic on each leaf of $\mathcal{F}$,
and the form $\sigma$ is non degenerate with $d'\sigma=0$; ii)
$\forall X\in\mathcal{V}^{1,0}_{pr}(M)$, the form $L_X\theta$
vanishes on $F$, iii) $\forall
X,Y\in\mathcal{V}^{1,0}_{pr}(M)$,one has
\begin{equation} \label{condHam} \flat_\theta(p_F[X,Y])=
d\{\sigma(X,Y)\} \end{equation} (i.e., $p_F[X,Y]$ is the
Hamiltonian vector field of the function
$\sigma(X,Y)=-\omega(X,Y)$).\end{corol} \noindent{\bf Proof.} Use
Proposition \ref{Vorobregular} with $G=\{0\}$. Q.e.d.
\begin{rem} \label{obscazsimplectic} {\rm
By looking at the bidegree of the various terms in the
decomposition of $d\omega$, we see that an $\mathcal{F}$-coupling
form $\omega$ given by (\ref{descsymplcoup}) is symplectic iff
\begin{equation} \label{eqsymplcoupl2}
d'\sigma=0,\;d''\sigma=-\partial\theta,\; d'\theta=0,\;
d''\theta=0.
\end{equation}

Furthermore, if we give the pair $(H,\theta_{0,2})$ where
$d'\theta=d''\theta=0$, $\partial\theta$ is a $d''$-closed
$(2,1)$-form, and $[\partial\theta] \in H^1(M,\Phi^2)$, where
$\Phi^2$ is the sheaf of germs of projectable cross sections of
$\wedge^2(ann\,F)$, is a cohomological obstruction to the
existence of a coupling form $\sigma_{2,0}$ associated with the
given pair such that the corresponding form $\omega$ of
(\ref{descsymplcoup}) is symplectic. (See \cite{V4} for the
cohomology space mentioned above.)}
\end{rem}
\par
Let $P''$ be a leaf-tangent bivector field on $(M,\mathcal{F})$. A
{\it transversally non degenerate extension} of $P''$ is a
bivector field $P$ such that $H=\sharp_P(ann\,F)$ is a normal
bundle of the foliation $\mathcal{F}$ and the corresponding
$(0,2)$-component of $P$ is the given bivector field $P''$. Then,
of course, $P$ is $\mathcal{F}$-coupling. An example is given by
\begin{prop} \label{existconexPoisson} If the foliation $F$ is a
locally trivial fibration over a $2q$-dimensional basis $B$ and if
$P''$ is a regular leaf-tangent Poisson structure on $(M,F)$,
there exists a transversally non degenerate extension $P$ of $P''$
which satisfies condition {\rm (\ref{condLie})}. \end{prop}
\noindent {\bf Proof.} By Weinstein's local structure theorem
(e.g., \cite{V1}), $M$ may be covered by open, connected,
neighborhoods endowed with local coordinates $(x^a, y^\nu,
p_\mu,q^\mu)$ such that $x^a=const.$ define $\mathcal{F}$,
$x^a=const.,y^\nu=const.$ define the symplectic foliation
$\mathcal{S}$ of $P''$, and $$P''=
\sum_\mu\frac{\partial}{\partial p_\mu}\wedge
\frac{\partial}{\partial q^\mu}.$$ Then,
$$\sum_{a=1}^{q}\frac{\partial}{\partial x^a}\wedge
\frac{\partial}{\partial
x^{a+q}}+ \sum_\mu\frac{\partial}{\partial p_\mu}\wedge
\frac{\partial}{\partial q^\mu}$$ are local
transversally non degenerate extensions of the required type. They
can be glued up to a global extension, which satisfies
(\ref{condLie}), by the pull-backs to $M$ of a partition of unity
of the basis $B$ of the fibration $\mathcal{F}$. Q.e.d.

In \cite{GLS}, the process started by Proposition
\ref{existconexPoisson} is continued up to the construction of a
coupling symplectic structure in the case where $P''$ is
leaf-tangent symplectic and the fibers are connected, simply
connected and compact. But, the difficulty  to get transversally
non degenerate, {\it Poisson extensions} of more general,
leaf-tangent, Poisson structures $P$ growth if the rank of $P$ is
lower. In particular, let us ask whether the structure $P''=0$ can
be extended to a coupling Poisson structure $P$. For $P''=0$,
(\ref{condLie}) is satisfied, and (\ref{condcurb}) shows that the
normal bundle $H$ must be integrable. Therefore extensions may
exist only on a locally product manifold with structural
foliations $\mathcal{F}$, $\mathcal{H}$, and, in view of
(\ref{d'sigmazero}), an extension is determined by a smooth family
$\sigma$ of symplectic structures of the leaves of $\mathcal{H}$.

It is natural to inquire about the existence of a coupling Poisson
structure (\ref{Pcoupling}) that is projectable onto the space of
leaves of the given foliation $\mathcal{F}$. Projectability holds
iff $(L_YP')(\alpha,\beta)=0$, $\forall Y\in\Gamma F$,
$\forall\alpha,\beta\in\Gamma(ann\,F)$, where, in fact, it is
enough to take projectable $1$-forms $\alpha,\beta$. The answer is
given by
\begin{prop} \label{P'foliat} $P$ given by {\rm
(\ref{Pcoupling})} is a projectable coupling Poisson bivector
field iff i) $P''$ is Poisson, ii) $H$ is integrable; iii) the
mod. $F$ equivalent $2$-form $\sigma$ of $P'$ is a transversal
symplectic form of $\mathcal{F}$; iv) {\rm (\ref{condLie})} holds.
\end{prop} \noindent{\bf Proof.} If $P'$ is projectable, so is the
equivalent $2$-form $\sigma$, and we see that (\ref{d'sigmazero})
implies iii) and that (\ref{condcurb}) implies the integrability
of $H$. Conversely, iii), iv) imply (\ref{d'sigmazero}),
(\ref{condLie}) and the projectability of $\sigma$ and $P'$,
which, together with ii), shows that the two sides of
(\ref{condcurb}) are zero. Q.e.d.

Thus, projectable, $\mathcal{F}$-coupling Poisson tensors exist
only on locally product manifolds $M$ with structural foliations
$\mathcal{F,H}$ where $\mathcal{H}$ is a leafwise symplectic
foliation.

\begin{rem} \label{tame} {\rm
If an almost coupling, projectable, bivector field $P$, given by
(\ref{Pcoupling}) is Poisson, $P'$ alone necessarily is a Poisson
bivector field, since by (\ref{descompSS1}) the first two
conditions (\ref{farabarP}) imply $[P',P']=0$. In the terminology
of \cite{V5}, $P'$ yields a {\it tame Hamiltonian structure} of
the foliation $\mathcal{F}$. Conversely, if we have such a tame
structure with a corresponding Poisson structure $P'$ of $M$, and
a leaf-tangent Poisson structure $P''$, and if the third condition
(\ref{farabarP}) holds, we get a projectable, almost coupling,
Poisson structure on $(M,\mathcal{F})$. Of course, we may always
take $P''=0$, hence, the tame Hamiltonian structures of
$\mathcal{F}$ and the projectable, almost coupling, Poisson
structures of $(M,\mathcal{F})$ are equivalent objects (not in a
one-to-one correspondence, however).} \end{rem}

\section{Vorobiev-Poisson structures}
In \cite{Vor}, Vorobiev extends an earlier construction of
Montgomery-Marsden-Ra\c tiu \cite{MMR} to a large class of bundles
of Lie coalgebras and obtains coupling Poisson structures in a
neighborhood of the zero section of these bundles. Here, we give a
simpler presentation of Vorobiev's results.

Let $p:\mathbb{G}^*\rightarrow B$ be a bundle of Lie coalgebras as
in Example
\ref{fibratePoisson} with the supplementary conditions: i) $B$ is
a symplectic manifold with the symplectic form $\omega$; ii) the
dual Lie algebras bundle $\mathbb{G}\rightarrow B$ is the kernel
of the (surjective) anchor $\rho:A\rightarrow TB$ of a transitive
Lie algebroid $p:A\rightarrow B$. (We refer the reader to
\cite{Mac} for the general theory of Lie algebroids.)

Let \begin{equation}\label{splittingA}
A=Q\oplus\mathbb{G}\end{equation} be a splitting of the vector
bundle $A$, and $p_Q,p_{\mathbb{G}}$ the corresponding natural
projections. Then $\rho|_Q:Q\rightarrow TB$ is an isomorphism and
we denote by $\gamma:TB\rightarrow Q$ its inverse. Obviously, we
have \begin{equation} \label{propluigamma} \gamma(\rho(s))=
p_Q(s), \hspace{1cm} \forall s\in A.\end{equation}

Since we will make some local coordinate checks later, we fix the
following notation now. Let $U\subseteq B$ be an open
local-trivialization neighborhood of all the vector bundles above
and $(x^i)$ $(i=1,...,n=dim\,B)$ local coordinates on $U$. Then,
we get the local basis ${\bf q}_i=\gamma(\partial/\partial x^i)$
of $Q$, and we may complete it by a local basis $({\bf g}_a)$ of
$\mathbb{G}$ $(a= 1,...,k=rank\,\mathbb{G})$  to a local basis of
$A$. We will denote by $(\theta^a)$ the dual basis of $({\bf
g}_a)$ for the vector bundle $\mathbb{G}^*$. These bases define
fiberwise local coordinates $(y^a)$ and $(y_a)$ on $\mathbb{G}$,
$\mathbb{G}^*$, respectively. Since $\mathbb{G}=ker\,\rho$ and
$\rho$ is a morphism of Lie algebroids, the local expressions of
the Lie bracket of $A$ must be of the form
\begin{equation} \label{LiebracketA} \begin{array}{l}
[{\bf g}_a,{\bf g}_b]_A=\alpha^c_{ab}(x){\bf g}_c,\hspace{2mm}
[{\bf g}_a,{\bf q}_i]_A=
\beta^c_{ai}(x){\bf g}_c,\vspace{2mm}\\

[{\bf q}_i,{\bf q}_j]_A=\gamma^c_{ij}(x){\bf g}_c +
\gamma_{ij}^h(x){\bf q}_h. \end{array} \end{equation} (Again, we
use the Einstein summation convention.) Finally, we indicate that
a free use of the following natural identifications will be made
below:
\begin{equation}\label{identificari}
\eta^a\frac{\partial}{\partial y^a}\Leftrightarrow \eta^a{\bf
g}_a,\;\;\mu_a\frac{\partial}{\partial y_a}\Leftrightarrow
\mu_a\theta_a. \end{equation}

Now, again since $\mathbb{G}=ker\,\rho$, the formula
\begin{equation} \label{conexpeM} \nabla_X\eta=[\gamma(X),\eta]_A,
\hspace{5mm}X\in\mathcal{V}^1(B),\;\eta\in\Gamma\mathbb{G},
\end{equation}
defines a connection of the vector bundle $\mathbb{G}$. In the
usual way, $\nabla$ yields a dual connection on $\mathbb{G}^*$ and
connections on all the associated tensor bundles of $\mathbb{G}$,
which we still denote by $\nabla$, except for situations where we
want to emphasize the connection on $\mathbb{G}^*$ and, then, we
will denote it by $\nabla^*$. Using (\ref{LiebracketA}), we see
that the local components of $\nabla$ are given by
\begin{equation} \label{coefdeconex}
\nabla_{\frac{\partial}{\partial x^i}}{\bf g}_a=\Gamma_{ai}^b{\bf
g}_b,\hspace{5mm}\Gamma_{ai}^b=-\beta_{ai}^b.\end{equation}

The Jacobi identity of the $A$-bracket is equivalent with
\begin{equation}
\label{Cparalel} \nabla C=0 \end{equation} for the tensor $C$
defined by (\ref{tensorC}). Furthermore, notice that
$$\rho([\gamma(X),\gamma (Y)]_A)=[X,Y]$$ implies
$$p_Q[\gamma(X),\gamma (Y)]_A=\gamma([X,Y]).$$ Then, by a
straightforward computation, we see that the curvature $R_\nabla$
satisfies the condition
\begin{equation} \label{curbinexemplu} R_\nabla(X,Y)\eta=
[p_{\mathbb{G}}[\gamma(X),\gamma(Y)]_A,\eta]_A.\end{equation}

Finally, we write
\begin{equation} \label{descompptnabla}
T\mathbb{G}^*=\mathcal{H}\oplus\mathcal{V}, \end{equation} with
the projections $p_{\mathcal{H}},p_{\mathcal{V}}$, where
$\mathcal{V}$ is tangent to the fibers and $\mathcal{H}$ is the
horizontal distribution of $\nabla^*$.

On $\mathbb{G}^*$, we have a triple $(P'',H,\sigma)$ as in
Proposition \ref{thVorobiev}, where $P''=\mathbb{L}$, $\mathbb{L}$
being the leaf-tangent Poisson bivector field defined by
(\ref{LiePinvariant}), $H=\mathcal{H}$ and
\begin{equation} \label{sigmapealgLie}
\sigma_z(\mathcal{X},\mathcal{Y})=\omega_{p(z)}(X,Y)-
z(p_{\mathbb{G}}[\gamma(X),\gamma(Y)]_A),\end{equation} where
$z\in\mathbb{G}^*$, $\mathcal{X},\mathcal{Y}\in\Gamma\mathcal{H}$
are the horizontal lifts of $X=p_*\mathcal{X},Y=p_*\mathcal{Y}$.
The definition of $\sigma$ is completed by asking it to have
bidegree $(2,0)$ with respect to (\ref{descompptnabla}).

The local expressions of the distribution $\mathcal{H}$ and the
form $\sigma$ follow from well known formulas of differential
geometry, and they are:
\begin{equation}\label{liftorizontal}\mathcal{H}=span
\{\mathcal{X}_i = \frac{\partial}{\partial x^i}
+\Gamma_{ai}^by_b\frac{\partial}{\partial y_a}\},\end{equation}
\begin{equation} \label{sigmalocal}
\sigma_z(\mathcal{X}_i,\mathcal{X}_j)=\omega_{ij}(x)-\gamma_{ij}^c
y_c,\end{equation} where $\mathcal{X}_i$ is the horizontal lift of
$\partial/\partial x^i$, the coefficients $\gamma,\Gamma$ are
those of (\ref{LiebracketA}) and (\ref{coefdeconex}), and
$\omega_{ij}$ are the natural local components of $\omega$ on $B$.
(Notice from (\ref{liftorizontal}) that the horizontal lifts of
vector fields of $B$ are projectable with respect to the vertical
foliation $\mathcal{V}$ of $\mathbb{G}^*$.)
\begin{prop}\label{structuraVor} {\rm\cite{Vor}}
There exists a neighborhood $\mathcal{U}$ of $B$, seen as the zero
section of $\mathbb{G}^*$, where the triple
$(\mathbb{L},\mathcal{H},\sigma)$ defines a coupling Poisson
bivector field.
\end{prop}
\noindent {\bf Proof.} Using (\ref{liftorizontal}), one gets
\begin{equation} \label{derivfunctliniare}
\mathcal{X}(<z,\eta>)=<z,\nabla_X\eta>\hspace{1cm}
(z\in\mathbb{G}^*, \eta\in\mathbb{G},\;p(z)=p(\eta)).
\end{equation} Accordingly,
$$d\sigma(\mathcal{X},\mathcal{Y},\mathcal{Z})=
\sum_{Cycl(\mathcal{X},\mathcal{Y},\mathcal{Z})}
\{\mathcal{X}\sigma(\mathcal{Y},\mathcal{Z})+
\sigma(\mathcal{X},[\mathcal{Y},\mathcal{Z}])\}$$
$$=d\omega(X,Y,Z)-\sum_{Cycl(X,Y,Z)}
<z,\nabla_X(p_{\mathbb{G}}[\gamma(Y),\gamma(Z)]_A)> $$ $$-
\sum_{Cycl(X,Y,Z)}
<z,p_{\mathbb{G}}[\gamma(X),p_Q[\gamma(Y),\gamma(Z)]_A]_A)>$$
$$=d\omega(X,Y,Z)-\sum_{Cycl(X,Y,Z)}
<z,p_{\mathbb{G}}[\gamma(X),[\gamma(Y),\gamma(Z)]_A)]_A> =0,$$
and, in view of the closedness of $\omega$ and of the Jacobi
identity for $A$, (\ref{d'sigmazero}) holds.

Furthermore, we have $$(L_{\mathcal{X}}\mathbb{L})_z
(\eta,\nu)=\mathcal{X}(<z,C(\eta,\nu>) -
\mathbb{L}_z(L_{\mathcal{X}}\eta,\nu) -
\mathbb{L}_z(\eta,L_{\mathcal{X}}\nu),$$ where $\eta,\nu\in\Gamma
T^*\mathbb{G}^*\Leftrightarrow \Gamma\mathbb{G}$. The Lie
derivatives in the right hand side may be evaluated at
$z\in\mathbb{G}^*$ by using the identification $z\Leftrightarrow
y_a(\partial/\partial y_a)$, and formulas (\ref{liftorizontal}),
(\ref{derivfunctliniare}):
$$<L_{\mathcal{X}}\eta,y_a\frac{\partial}{\partial y_a}>=\mathcal{X}<z,\eta>
-<\eta,[\mathcal{X},y_a\frac{\partial}{\partial y_a}]>=
<z,\nabla_X\eta>.$$

As a consequence, we obtain
\begin{equation} \label{derivorizC} (L_{\mathcal{X}}\mathbb{L})_z
(\eta,\nu)=<z,(\nabla_XC)(\eta,\nu)>,\end{equation} and
(\ref{condLie}) follows from (\ref{Cparalel}).

Finally, (\ref{condcurb}) can be checked as follows. From
(\ref{PpemathcalG}), (\ref{tensorC}), and (\ref{curbinexemplu}) we
get
$$\mathbb{L}_z(p_{\mathbb{G}}[\gamma(X),\gamma(Y)]_A,\eta)
=<z,[p_{\mathbb{G}}[\gamma(X),\gamma(Y)]_A,\eta]_A>$$
$$=<z,R_\nabla(X,Y)\eta>=-<R_{\nabla^*}(X,Y)z,\eta>,$$
where again $z\Leftrightarrow y_a(\partial/\partial y_a)$.
Equivalently,
$$\sharp_{\mathbb{L}_z}(p_{\mathbb{G}}[\gamma(X),\gamma (Y)]_A)
=-R_{\nabla^*}(X,Y)z.$$

But, by classical connection theory, $-R_{\nabla^*}(X,Y)z$ is
equal to $p_\mathcal{V}[\mathcal{X},\mathcal{Y}](z)$, here
identified with a vector of the fiber $\mathbb{G}_{p(z)}$. On the
other hand, $\forall \mu\in\Gamma\mathbb{G}$ the $(0,1)$-component
of $d<z,\mu_z>$ with respect to (\ref{descompptnabla}) identifies
with $\mu$. Accordingly, (\ref{sigmapealgLie}) yields
$$\sharp_{\mathbb{L}_z}
\{d[\sigma(\mathcal{X},\mathcal{Y})]\}=P_{\mathcal{V}}
[\mathcal{X},\mathcal{Y}], $$ which is (\ref{condcurb}) in our
case.

Therefore, since $\sigma$ is non degenerate on a neighborhood
$\mathcal{U}$ of the zero section of $\mathbb{G}^*$, there exists
a corresponding coupling Poisson structure on $\mathcal{U}$.
Q.e.d.
\begin{defin}\label{defstrVorobiev} {\rm The coupling
Poisson structure defined by Proposition \ref{structuraVor} on the
neighborhood $\mathcal{U}$ of $B$ will be called a {\it
Vorobiev-Poisson structure}.}\end{defin}
\begin{rem}\label{obscoiso} {\rm \cite{Vor}
If the distribution $$\mathcal{D}= \{X\in
TB\,/\,p_{\mathbb{G}}[\gamma(X),\gamma(Y)]_A=0,\, \forall Y\in
TB\}$$ is coisotropic with respect to $\omega$, the
Vorobiev-Poisson structure is globally coupling (i.e.,
$\mathcal{U}=\mathbb{G}^*$). Indeed, the local equations of
$\mathcal{D}$ are \begin{equation}
\label{distribD}\gamma_{ij}^aX^i=0\hspace{1cm}(X=X^i\frac{\partial}{\partial
x^i}) ,\end{equation} and $\mathcal{D}$ is coisotropic iff the
equations (\ref{distribD}) together with $\omega_{ij}X^iY^j=0$
imply $\gamma_{ij}^aY^i=0$. Now, assume that
$\mathcal{X}=X^i\mathcal{X}_i\in ann\,\sigma_z$ i.e.,
\begin{equation} \label{anihilsigma}
(\omega_{ij}-y_c\gamma^c_{ij})X^iY^j=0,\hspace{5mm}
\forall\mathcal{Y}=Y^j\mathcal{X}_j.\end{equation} Then, we get
$\omega_{ij}X^iY^j=0$ for all the vectors $Y\in\mathcal{D}_{p(z)}$
and $X=X^i(\partial/\partial x^i)$ must belong to $\mathcal{D}$.
This reduces (\ref{anihilsigma}) to $X\in ann\,\omega_{p(z)}$ and,
since $\omega$ is non degenerate, we get $X=0$ and
$\mathcal{X}=0$, i.e., $\sigma_z$ is non degenerate. }\end{rem}

Furthermore, one has the following important result.
\begin{prop}\label{echivalenta}
{\rm\cite{Vor}} The Vorobiev-Poisson structures defined by two
splittings \begin{equation}\label{twosplittings}
A=Q\oplus\mathbb{G},\;\;A=\tilde Q\oplus\mathbb{G} \end{equation}
on neighborhoods $\mathcal{U}_1,\mathcal{U}_2$ of $B$ in
$\mathbb{G}^*$ are Poisson-equivalent in a neighborhood
$V\subseteq \mathcal{U}_1\cap \mathcal{U}_2$.
\end{prop} \noindent {\bf Proof.} The notation below is that of
Proposition \ref{structuraVor} with the addition of a tilde for
everything related to the second splitting .

The difference $\phi=\tilde\gamma-\gamma$ is a $\mathbb{G}$-valued
$1$-form on $B$ and, if $s\in\Gamma A$ has the two decompositions
$$s=p_{\mathbb{G}}(s)+p_{Q}(s)
=\tilde p_{\mathbb{G}}(s)+p_{\tilde Q}(s),$$ (\ref{propluigamma})
implies
\begin{equation}\label{legatura12}
\tilde p_{\mathbb{G}}(s)=p_{\mathbb{G}}(s)
-\phi(\rho(s)),\;\;p_{\tilde
Q}(s)=p_Q(s)+\phi(\rho(s)).\end{equation} Furthermore, on
$\mathbb{G}^*$ we get a scalar $1$-form $\psi\in ann\,\mathcal{V}$
by means of the formula
\begin{equation}\label{formapsi}
\psi_z(\mathcal{X})=<z,\phi(X)> ,\hspace{5mm}
z\in\mathbb{G}^*.\end{equation}

Using $\phi$, we may define a homotopy, i.e., a family of
splittings $A=Q_t\oplus\mathbb{G}$, given by the projectors
\begin{equation}\label{homotopiadescomp}p_{Q_t} =p_Q+t(\phi\circ\rho),
\hspace{1cm}t\in\mathbb{R},\end{equation} which is such that
$Q_0=Q$, $Q_1=\tilde Q$, with the corresponding homotopy class
$P_t$ of Vorobiev-Poisson bivector fields defined by the triples
$(\mathbb{L},\mathcal{H}_t,\sigma_t)$ of the connection $\nabla^t$
associated to (\ref{homotopiadescomp}) via (\ref{conexpeM}).

With (\ref{coefdeconex}) and the definition of $\phi$, it is easy
to compute the connection coefficients of $\nabla^t$ and get
\begin{equation}\label{liftorizontalt}\mathcal{H}_t=span
\{\mathcal{X}_{t,i} = \mathcal{X}_{i}
+t\alpha_{ac}^b\phi^c_iy_b\frac{\partial}{\partial
y_a}\},\end{equation} where the components $\phi^c_i$ are given by
\begin{equation}\label{compluiphi}
\phi(\frac{\partial}{\partial x^i})=\phi^c_i{\bf
g}_c.\end{equation} The corresponding basis of
$ann\,\mathcal{H}_t=\mathcal{V}^*$ consists of the forms
\begin{equation}\label{cobazat} \mu_{t,a}=\mu_a
-t\alpha_{ac}^b\phi_i^c y_bdx^i,\end{equation} where
$\mu_a=\mu_{0,a}$. From (\ref{liftorizontalt}) we get the
horizontal lift of $X\in TB$ to $\mathcal{H}_t$
\begin{equation}\label{relintreliftorizontal} \mathcal{X}_t(z)=
\mathcal{X}(z)-t\,coad_{\phi(X)}(\mathbb{E}(z)),\hspace{3mm}
z\in\mathbb{G}^*,\, \mathbb{E}(z)=y_a\frac{\partial}{\partial
y_a}\Leftrightarrow y_a\theta^a \end{equation} ($\mathbb{E}$ is
the infinitesimal homothety of $\mathbb{G}^*$). Notice also that
(\ref{cobazat}) may be seen as a bijection
$\lambda\mapsto\lambda_t$ between $ann\,\mathcal{H}$ and
$ann\,\mathcal{H}_t$ given by:
\begin{equation} \label{corespaanihilatorilor}
\lambda_t=\lambda-tL_{\sharp_{\mathbb{L}}\lambda}\psi,
\end{equation}
where $\psi$ is the $1$-form (\ref{formapsi}).

Now, we can compute the $2$-form $\sigma_t$ of the triple
$(\mathbb{L},\mathcal{H}_t,\sigma_t)$. Since the horizontal lifts
are such that the differences $\tilde{\mathcal{X}}-\mathcal{X},
\tilde{\mathcal{Y}}-\mathcal{Y}$ are vertical, and using
(\ref{legatura12}), we get
\begin{equation}\label{relintresigma}
\sigma_{t}(z)(\mathcal{X}_t,\mathcal{Y}_t) =
\sigma_{t}(z)(\mathcal{X},\mathcal{Y}) =
\sigma(z)(\mathcal{X},\mathcal{Y})
\end{equation}
$$ - <z,t[\gamma(X),\phi(Y)]_A- t[\gamma(Y),\phi(X)]_A -
t\phi([X,Y])+ t^2[\phi(X),\phi(Y)]_A>.$$

Furthermore, using (\ref{derivfunctliniare}),
(\ref{relintresigma}) becomes
\begin{equation}\label{relintresigma2}
\sigma_{t}(\mathcal{X}_t,\mathcal{Y}_t) =
\sigma(\mathcal{X},\mathcal{Y})-td\psi(\mathcal{X},\mathcal{Y})
-t^2\mathbb{L}(\phi(X),\phi(Y)).
\end{equation}

As in \cite{Vor}, at this point we define a time-dependent vector
field $\Xi_t\in\Gamma(\mathcal{H}_t)$ by
\begin{equation}\label{campXit} \Xi_t=\sharp_{P_t}\psi.
\end{equation} Equivalently, and with the notation of
(\ref{Pcoupling}), we have \begin{equation} \label{campXit2}
\flat_{\sigma_t}\Xi_t=i(\Xi_t)\sigma_t=-\psi.\end{equation} The
vector field $\Xi_t$ is equivalent with the autonomous vector
field $\tilde{\Xi}=\Xi_t+\partial/\partial t$ on
$\mathbb{G}^*\times\mathbb{R}$, and we will prove that the flow
$\Phi_t$ of $\tilde{\Xi}$ preserves the lift of the tensor field
$P_t$ to $\mathbb{G}^*\times\mathbb{R}$. Then the projection of
the diffeomorphism $\Phi_1$ onto $\mathbb{G}^*$ will be the
required equivalence of coupling Poisson structures.

Thus, the proof will be accomplished if we show that
\begin{equation} \label{LiePtzero} L_{\tilde{\Xi}}P_t=0.
\end{equation} Condition (\ref{LiePtzero}) obviously holds if one
of the arguments is $dt$. In the other cases, we shall see that
(\ref{LiePtzero}) is a consequence of the fact that $P_t$
satisfies the characteristic conditions (\ref{farabarP}) of the
coupling Poisson tensors for all $t$.

Indeed, the third condition (\ref{farabarP}) for $\alpha=\psi$
yields $$(L_{\Xi_t}P_t) (\lambda_t,\mu_t)=0$$ and
(\ref{corespaanihilatorilor}) yields
$$L_{\frac{\partial}{\partial t}}P_t(\lambda_t,\mu_t)=
\frac{\partial}{\partial t}(\mathbb{L}(\lambda,\mu)=0.$$ Therefore
(\ref{LiePtzero}) holds for arguments in $ann\,\mathcal{H}_t$.

For mixed arguments $\alpha\in ann\,\mathcal{V},\lambda_t\in
ann\,\mathcal{H}_t$ we get
$$(L_{\Xi_t}P_t)(\alpha,\lambda_t)=
\alpha([\sharp_{P''_t}\lambda_t,\sharp_{P'_t}\psi])
-\lambda_t([\sharp_{P'_t}\alpha,\sharp_{P'_t}\psi])$$ and, if the
second term is replaced in agreement with the second condition
(\ref{farabarP}) and $\lambda_t$ is expressed by
(\ref{corespaanihilatorilor}), the result is
$$(L_{\Xi_t}P_t)(\alpha,\lambda_t)=
P'_t(L_{\sharp_{P''_t}\lambda_t}\psi,\alpha)
=P'_t(L_{\sharp_{\mathbb{L}}\lambda}\psi,\alpha).$$ On the other
hand, using (\ref{corespaanihilatorilor}) again, we get
$$(L_{\frac{\partial}{\partial t}}P_t)(\alpha,\lambda_t)=
-P_t(\alpha,\frac{\partial\lambda_t}{\partial t})=
P'_t(\alpha,L_{\sharp_{\mathbb{L}}\lambda}\psi).$$ Therefore,
(\ref{LiePtzero}) holds for mixed arguments.

Finally, for $\alpha,\beta\in ann\,\mathcal{V}$, and if
$\mathcal{X}_t=\sharp_{P'_t}\alpha,\mathcal{Y}_t=\sharp_{P'_t}\beta$,
the first condition (\ref{farabarP}) for $\gamma=\psi$ yields
$$(L_{\Xi_t}P_t)(\alpha,\beta)=
d\psi(\mathcal{X}_t,\mathcal{Y}_t).$$ Furthermore, with
(\ref{liftorizontalt}), we deduce
$$(L_{\Xi_t}P_t)(\alpha,\beta)=
d\psi(\mathcal{X},\mathcal{Y})+2t\mathbb{L}(\phi(X),\phi(Y)),$$
where the notation is that of the formulas
(\ref{relintreliftorizontal}), (\ref{relintresigma2}). Then,
$$(L_{{\partial}{\partial t}}P_t)(\alpha,\beta)=
\frac{\partial}{\partial t}[P_t(\alpha,\beta)]=
\frac{\partial}{\partial
t}[\sigma_t(\mathcal{X}_t,\mathcal{Y}_t)]$$
$$=-d\psi(\mathcal{X},\mathcal{Y})-2t\mathbb{L}(\phi(X),\phi(Y)).$$
Therefore, (\ref{LiePtzero}) also holds for two arguments in
$ann\,\mathcal{V}$. Q.e.d.
\begin{rem} \label{adouastrVorob} {\rm
It is important to notice that the equivalence of Poisson
structures given by Proposition \ref{echivalenta} does not
preserve the foliation $\mathcal{V}$. On the other hand, we also
notice that the triple $(-\mathbb{L},H, \tilde\sigma)$, where
\begin{equation} \label{tildasigmapealgLie}
\tilde{\sigma}_z(\mathcal{X},\mathcal{Y})=\omega_{p(z)}(X,Y)+
z(p_{\mathbb{G}}[\gamma(X),\gamma(Y)]_A),\end{equation} is also a
coupling Poisson structure on a neighborhood of $B$ in
$\mathbb{G}^*$. This follows from Remark \ref{obssuplementara}
since $-\tilde\sigma=\sigma-2p^*\omega$. The new structure also
satisfies the equivalence property of Proposition
\ref{echivalenta}. The proof is the same except for a change of
sign in the definition of the vector field $\Xi_t$. This structure
will also be called a Vorobiev-Poisson structure.}\end{rem}

Following \cite{Vor}, it is possible to apply the previous
construction to an embedded leaf $S$ of a Poisson manifold
$(M,P)$. The Lie algebroid structure of the cotangent bundle
$T^*M$ restricts to a transitive Lie algebroid
$T^*M|_S\longrightarrow S$ and the kernel of the anchor of this
algebroid is the {\it conormal bundle} of $S$, i.e., the
annihilator of $TS$ in $T^*M|_S$. Let $NS$ be a normal bundle of
$S$ (i.e., $TM|_S=TS\oplus NS$) and $U$ a tubular neighborhood of
$S$ with the fibers tangent to $NS$. At the points of $S$ there
exist local adapted coordinates $(x^\alpha,x^\kappa)$
$(\alpha=1,...,codim(S);\kappa=codim(S)+1,...,dim(M))$ such that
the local equations of $S$ are $x^\alpha=0$ and
\begin{equation}\label{bazeinlungfoaie}
N^*S=ann(TS)= span\{dx^\alpha|_{S}\},\;\;T^*S=ann(NS)=
span\{dx^\kappa|_{S}\}.\end{equation} The vector bundles $N^*S$
and $T^*S$ may play the role of $\mathbb{G}$ and $Q$ of Vorobiev's
construction, respectively, and the bases (\ref{bazeinlungfoaie})
may have the role of the bases $({\bf g}_a)$, $({\bf q}_i)$ of
(\ref{LiebracketA}). If
$P^{\alpha\beta},P^{\alpha\kappa},P^{\kappa\nu}$ are the local
components of $P$ with respect to the local coordinates defined
above, one has $P^{\alpha\beta}|_S=P^{\alpha\nu}|_S=0$ , whence
$$\frac{\partial P^{\alpha\beta}}{\partial x^\kappa}|_{S}=0,\;\;
\frac{\partial P^{\alpha\nu}}{\partial x^\kappa}|_{S}=0.$$ Accordingly, the
brackets (\ref{LiebracketA}) of the present case will be
\begin{equation}\label{croseteptaproxlin} \begin{array}{l}
\{dx^\alpha|_{S},dx^\beta|_{S}\}=\frac{\partial
P^{\alpha\beta}}{\partial x^\gamma}|_{S}dx^\gamma|_{S},
\,\,\{dx^\alpha|_{S},dx^\kappa|_{S}\}=\frac{\partial
P^{\alpha\kappa}}{\partial x^\gamma}|_{S}dx^\gamma|_{S},
\vspace{2mm}\\
\{dx^\kappa|_{S},dx^\nu|_{S}\}= \frac{\partial
P^{\kappa\nu}}{\partial x^\gamma}|_{S}dx^\gamma|_{S} +
\frac{\partial P^{\kappa\nu}}{\partial
x^\theta}|_{S}dx^\theta|_{S}.\end{array} \end{equation} (In the
formulas above, $\beta,\gamma$ have the same domain as $\alpha$
and $\nu,\theta$ have the same domain as $\kappa$.)

Therefore, one has a Vorobiev-Poisson structure on a neighborhood
of $S$ and the local formulas (\ref{PpemathcalG}),
(\ref{liftorizontal}) and (\ref{sigmalocal}) show that the
associated triple is given by
\begin{equation}\label{liniarizare} \begin{array}{l}
P''=\frac{1}{2}\xi^\gamma\frac{\partial P^{\alpha\beta}}{\partial
x^\gamma}|_{x^\gamma=0}\frac{\partial}{\partial\xi^\alpha}
\wedge\frac{\partial}{\partial\xi^\beta},\vspace{2mm}\\
H=span\left\{\frac{\partial}{\partial x^\kappa} -
\xi^\gamma\frac{\partial P^{\alpha\kappa}}{\partial
x^\gamma}|_{x^\gamma=0}\frac{\partial}{\partial\xi^\alpha}
\right\},\vspace{2mm}\\
\sigma=\frac{1}{2}(p_{\kappa\nu}- \xi^\gamma\frac{\partial
P^{\kappa\nu}}{\partial x^\gamma}|_{x^\gamma=0}) dx^\kappa\wedge
dx^\nu,
\end{array} \end{equation} where $\xi^\alpha$ are fiber coordinates
in the normal bundle $NS$,
$p_{\kappa\nu}P^{\nu\theta}=\delta_\kappa^\theta$ and the indices
take the same values as in (\ref{croseteptaproxlin}). The
invariant expression of the $2$-form $\sigma$ of
(\ref{liniarizare}) is
\begin{equation}
\label{approxlinPs}
\sigma_Z(\mathcal{X},\mathcal{Y})=\omega_{p(Z)}(X,Y)-
<Z,\{\xi,\eta\}_{P}>, \end{equation} where $Z\in NS$, $\xi,\eta\in
T^*S$, $X=\sharp_P\xi,Y=\sharp_P\eta$ and the bracket of $1$-forms
is that given by (\ref{bracket1}).

The Vorobiev-Poisson structure (\ref{liniarizare}) is defined up
to Poisson equivalence and its $S$-transversal part may be seen as
a linear approximation of the $S$-transversal part of the original
Poisson structure $P$.
\section{Jacobi structures on foliated manifolds}
In this section we discuss similar problems for Jacobi structures.
We recall that a Jacobi structure is a Lie bracket of the local
type on the algebra of differentiable functions on a manifold,
$C^\infty(M)$. Such a bracket must be of the form
\begin{equation}\label{Jacobibr} \{f,g\}=\Lambda(df,dg)+f(Eg)-g(Ef)
\hspace{5mm}(f,g\in C^\infty(M))\end{equation} where $\Lambda
\in\mathcal{V}^2(M),E\in\mathcal{V}^1(M)$ satisfy the conditions
\begin{equation}\label{LLzero}
[\Lambda,\Lambda]=2E\wedge\Lambda,\;\;L_E\Lambda=0.\end{equation}
Accordingly, we refer to such a pair $(\Lambda,E)$ as a {\it
Jacobi structure}. For a Jacobi structure,
$im\,\sharp_\Lambda+span\{E\}$ is an integrable generalized
distribution, called the {\it characteristic distribution}, with
either locally conformal symplectic or contact leaves. As general
references on Jacobi structures, we quote \cite{DLM} and
\cite{V3}.

We also recall the fundamental fact that the Jacobi structure
$(\Lambda,E)$ on the manifold $M$ is equivalent with the {\it
homogeneous} Poisson structure
\begin{equation} \label{homogPsptJacobi} P=e^{-t}\left(
\Lambda+\frac{\partial}{\partial t}\wedge E\right),\hspace{5mm}
t\in\mathbb{R},\end{equation} on $M\times\mathbb{R}$
\cite{{DLM},{V3}}.

The general conditions for a couple $(\Lambda,E)$ to define a
Jacobi structure on the foliated manifold $(M,\mathcal{F})$ with a
normal bundle $H$ of $\mathcal{F}$ may be written down using
(\ref{descompSS1})-(\ref{descompSS6}) to express (\ref{LLzero}).
\begin{defin} \label{precoupling} {\rm A pair $(\Lambda,E)$
that consists of a bivector field $\Lambda$ and a vector field $E$
on a manifold $M$ is {\it pre-coupling} with respect to a
foliation $\mathcal{F}$ if $\sharp_\Lambda(ann\,F)\cap F=\{0\}$
($F=T\mathcal{F}$). If $E$ is tangent to $\mathcal{F}$ and
$\Lambda\in\Gamma\wedge^2F$ the pair $(\Lambda,E)$ is said to be a
{\it leaf-tangent pair}.}\end{defin}

In the remaining part of the paper all the pairs $(\Lambda,E)$ on
$(M,\mathcal{F})$ will be pre-coupling. In particular, a
leaf-tangent pair obviously is pre-coupling.
\begin{defin}\label{cazuriJacobi} {\rm 1) The pre-coupling pair
$(\Lambda,E)$ is {\it of the $\mathcal{F}$-tangent} or {\it the
$\mathcal{F}$-normal type} if the vector field $E$ is always
tangent, respectively, never tangent to $\mathcal{F}$. 2)
$(\Lambda,E)$ is $\mathcal{F}$-{\it almost coupling via} $H$ if
$\sharp_\Lambda(ann\,F)\subseteq H$ for the normal bundle $H$ of
$\mathcal{F}$.  3) $(\Lambda,E)$ is {\it $\mathcal{F}$-coupling of
the first kind} if it is of the tangent type and
$\sharp_\Lambda(ann\,F)$ is a normal bundle $H$ of $\mathcal{F}$.
4) $(\Lambda,E)$ is {\it $\mathcal{F}$-coupling of the second
kind} if $\sharp_\Lambda(ann\,F)$ is a normal bundle $H$ of
$\mathcal{F}$ that contains $E$. 5) $(\Lambda,E)$ is {\it
$\mathcal{F}$-coupling of the third kind} if it is of the normal
type and $H'=\sharp_\Lambda[ann(F\oplus span\{E\})]$ is a
complementary distribution of $F\oplus(span\{E\})$.}
\end{defin}

In the first and second coupling cases, we always take the normal
bundle $H$ defined by the coupling condition and in the third case
we always take $H=H'\oplus span\{E\}$. For any foliation
$\mathcal{F}$ and any bivector field $\Lambda$ we have
$$F\cap\sharp_\Lambda(ann\,F)=\sharp_\Lambda(ker\,\Lambda|_{ann\,F}).
$$ Hence, $\sharp_\Lambda(ann\,F)$ is normal to $\mathcal{F}$ iff
\begin{equation}\label{eq1rang}
rank(\Lambda|_{ann\,\mathcal{F}})=codim\mathcal{F}.\end{equation}
Thus, if $E\in\Gamma F$, condition (\ref{eq1rang}) is equivalent
with coupling of the first kind while, on the contrary, if $E\in
im\,\sharp_\Lambda$, (\ref{eq1rang}) is equivalent with second
kind coupling. Similarly, we see that the third kind coupling
condition holds iff $E$ is never tangent to $\mathcal{F}$ and
\begin{equation}\label{eq2rang}
rank(\Lambda|_{ann(\mathcal{F}\oplus
span\{E\})})=codim\mathcal{F}-1.\end{equation}
\begin{example} \label{exfoiinJacobi} {\rm Assume that a
{\it characteristic leaf} $L$ of the Jacobi manifold
$(M,\Lambda,E)$ is transversal to the foliation $\mathcal{F}$. If
$L$ is locally conformal symplectic, the condition that $E$ is not
in $F$ and (\ref{eq1rang}) hold along $L$ hence, by continuity,
also hold in a neighborhood $U$ of $L$. Therefore, $(\Lambda,E)$
is $\mathcal{F}$-coupling of the second kind on $U$. If $L$ is a
contact manifold, the conditions for $(\Lambda,E)$ to be coupling
of the third kind hold along $L$ hence, also on a neighborhood $U$
of $L$. In particular, for any Jacobi structure, any embedded
even-dimensional leaf $L$, has a tubular neighborhood $U$ where
the Jacobi structure is coupling of the second kind with the
fibers of the tubular structure and any embedded odd-dimensional
leaf $L$ has a tubular neighborhood $U$ of $L$ where the Jacobi
structure is coupling of the third kind with the fibers of the
tubular structure.}\end{example}
\begin{example} \label{exIglesias} {\rm Let $\mathbb{G}\longrightarrow
B$ be a bundle of Lie algebras, which is the kernel of a
transitive Lie algebroid over a symplectic manifold, and assume
that there exists a global cross section $\zeta$ of the dual
bundle $\mathbb{G}^*$ that vanishes on the derived algebras of the
fibers of $\mathbb{G}$. Then $\zeta$ may be seen as a vertical
vector field on $\mathbb{G}^*$, and straightforward computations
show that the pair
\begin{equation}
\label{eqexIglesias}
\Lambda=\mathbb{L}+\mathbb{E}\wedge\zeta,\;\;E=-\zeta,\end{equation}
where $\mathbb{L},\mathbb{E}$ are defined by
(\ref{LiePinvariant}),  (\ref{relintreliftorizontal}), is a
leaf-tangent Jacobi structure on $\mathbb{G}^*$ with the vertical
foliation $\mathcal{V})$. (This is a particular case of the
general result of \cite{MIgl}.) Furthermore, if $P$ is a
Vorobiev-Poisson structure on a neighborhood of $B$ in
$\mathbb{G}^*$ with the vertical part $\mathbb{L}$, the pair
$\Lambda_P=P+\mathbb{E}\wedge\zeta,\;\;E=-\zeta$ is coupling of
the first kind on the same neighborhood, which is not a Jacobi
structure, however.} \end{example}

It is clear that coupling implies almost coupling, and coupling of
either the second or the third kind implies the normal type
condition. Coupling of the second kind may exist only if the
codimension of $\mathcal{F}$ is even and coupling of the third
kind may exist only if the codimension of $\mathcal{F}$ is odd.
Obviously, condition (\ref{eq1rang}) holds iff the bivector field
$P$ defined on $M\times\mathbb{R}$ by (\ref{homogPsptJacobi}) is
coupling for the foliation $\mathcal{F}\times \mathbb{R}$.
Therefore, this property of $P$ holds for $(\Lambda,E)$ coupling
of the first or second kind, and the two situations are
differentiated by $E\in\Gamma F$, $E\in
\sharp_\Lambda(ann\,F)$, respectively. Finally, if $(\Lambda,E)$
is of the normal type, coupling of the third kind holds iff $P$ is
coupling for the pullback of $\mathcal{F}$ to $M\times\mathbb{R}$
and, then, the corresponding normal bundle is $\tilde H\oplus
span\{\frac{\partial}{\partial t}\}$, where $\tilde H$ is the
pullback of $\sharp_\Lambda[ann(F\oplus span\{E\})]$ to
$M\times\mathbb{R}$.

Furthermore, we remark that a conformal change
\begin{equation}\{f,g\}^a=\frac{1}{a}
\{af,ag\}\hspace{5mm}(a=e^b,\;b\in C^\infty(M))\end{equation}
of an $\mathcal{F}$-coupling Jacobi structure of either the second
or the third kind $\{f,g\}$ leads to a Jacobi structure of the
same type. Indeed, if the original structure is given by
$(\Lambda,E)$, the new structure is given by \cite{{DLM},{V3}}
\begin{equation}\label{conformJ}
\Lambda^a=a\Lambda,\;E^a=aE+\sharp_\Lambda(da).
\end{equation} If the original structure is coupling of the first kind,
(\ref{conformJ}) is coupling of the first kind iff $a=const.$

With the definitions above, we get \begin{prop} \label{Jpefol1} On
a foliated manifold $(M,\mathcal{F})$, a pair $(\Lambda,E)$ of the
tangent type, which is almost coupling via a normal bundle $H$, is
a Jacobi structure iff
\begin{equation}\label{Lbarzero}
\Lambda=\Lambda'_{2,0}+\Lambda''_{0,2}\end{equation} with respect to
the bigrading defined by $H$ and
\begin{equation} \label{farabarJac} \begin{array}{l}
 d'\gamma(\sharp_{\Lambda'}\alpha,\sharp_{\Lambda'}\beta)
-(L_{\sharp_{\Lambda'}\gamma}\Lambda')(\alpha,\beta)=0,\vspace{2mm}\\

(L_{\sharp_{\Lambda''}\lambda}\Lambda')
(\alpha,\beta)+\lambda([\sharp_{\Lambda'}\alpha,
\sharp_{\Lambda'}\beta])=
\lambda(E)\Lambda'(\alpha,\beta),\vspace{2mm}\\

(L_{\sharp_{\Lambda'}\alpha}\Lambda'')(\lambda,\mu)= 0,
\vspace{2mm}\\

 d''\nu(\sharp_{\Lambda''}\lambda,\sharp_{\Lambda''}\mu)
-(L_{\sharp_{\Lambda''}\nu}\Lambda'')(\lambda,\mu)=
(E\wedge\Lambda'')(\lambda,\mu,\nu)
\end{array} \end{equation}
\begin{equation}\label{LieinJac1} L_{E}\Lambda''=0,\;
L_{E}\Lambda'=0\end{equation} for $1$-forms $\alpha,\beta,\gamma$
of bidegree $(1,0)$ and $\lambda,\mu,\nu$ of bidegree $(0,1)$.
\end{prop} \noindent{\bf Proof.} From
(\ref{descompSS1})-(\ref{descompSS6}), and since $E\in\Gamma F$,
we see that the first condition (\ref{LLzero}) is equivalent to
(\ref{Lbarzero}). Then, by looking at $L_E\Lambda=0$ for arguments
of bidegree $(1,0),(0,1)$ we get the first condition
(\ref{LieinJac1}) first, and then the second. Q.e.d.

Together, the last condition (\ref{farabarJac}) and the first
condition (\ref{LieinJac1}) are equivalent with the fact that
$(\Lambda'',E)$ is a leaf-tangent Jacobi structure on
$(M,\mathcal{F})$, and a pair $(\Lambda_{0,2},E_{0,1})$ is a
leaf-tangent Jacobi structure iff these two conditions hold.

Furthermore, if $(\Lambda=\Lambda'+\Lambda'',E)$ is coupling of
the first kind, we may use the $(2,0)$-form $\sigma$ that is
equivalent, mod. $F$ with $\Lambda'$, as we did in the Poisson
case (see Proposition \ref{thVorobiev}), and we get
\begin{prop} \label{caz1Jcusigma} If the pair
$(\Lambda,E)$ is coupling of the first kind, it is a Jacobi
structure iff $(\Lambda'',E)$ is leaf-tangent Jacobi and
\begin{equation} \label{VorobJac1} \begin{array}{c}
d'\sigma=0,\,L_X\Lambda''=0, \vspace{2mm}\\
-p_F[X,Y]=\sharp_{\Lambda''}\{d[\sigma(X,Y)]\}-\sigma(X,Y)E,\vspace{2mm}\\
L_E\Lambda'=0.
\end{array}\end{equation} for all the vector fields
$X,Y\in\mathcal{V}^{1,0}_{pr}(M).$\end{prop}

Notice also that the last condition (\ref{VorobJac1}) is
equivalent with \begin{equation} \label{ultimaJac1}
E(\sigma(X,Y))=0,\;[E,X]=0\hspace{5mm}
X,Y\in\mathcal{V}^{1,0}_{pr}(M).\end{equation}

Similarly, we get   \begin{prop} \label{Jpefol2} An almost
coupling pair $(\Lambda,E)$ of the normal type is a Jacobi
structure iff {\rm(\ref{Lbarzero})} holds and the following
conditions are satisfied:
\begin{equation} \label{farabarJac2} \begin{array}{l}
 d'\gamma(\sharp_{\Lambda'}\alpha,\sharp_{\Lambda'}\beta)
-(L_{\sharp_{\Lambda'}\gamma}\Lambda')(\alpha,\beta)=
E\wedge\Lambda'(\alpha,\beta,\gamma),\vspace{2mm}\\

(L_{\sharp_{\Lambda''}\lambda}\Lambda')
(\alpha,\beta)+\lambda([\sharp_{\Lambda'}\alpha,
\sharp_{\Lambda'}\beta])=0,
\vspace{2mm}\\

(L_{\sharp_{\Lambda'}\alpha}\Lambda'')(\lambda,\mu)=
\alpha(E)\Lambda''(\lambda,\mu),
\vspace{2mm}\\

 d''\nu(\sharp_{\Lambda''}\lambda,\sharp_{\Lambda''}\mu)
-(L_{\sharp_{\Lambda''}\nu}\Lambda'')(\lambda,\mu)=0,
\end{array} \end{equation}
\begin{equation}\label{LieinJac2}
(L_{E}\Lambda')(\alpha,\beta)=0,\;(L_{E}\Lambda'')(\lambda,\mu)=0,\;
(L_{E}\Lambda')(\alpha,\lambda)=(L_{E}\Lambda'')(\lambda,\alpha).
\end{equation}
\end{prop}

In the situation of Proposition \ref{Jpefol2}, the last condition
(\ref{farabarJac2}) is equivalent with the fact that $\Lambda''$
is a leaf-tangent Poisson structure on $(M,\mathcal{F})$ and $E$
is an infinitesimal automorphism of the former. If $E$ is a
projectable vector field conditions (\ref{LieinJac2}) become again
(\ref{LieinJac1}). Furthermore, in the case of a coupling pair of
the second kind, we may again use a $(2,0)$-form $\sigma$ mod.
$F$-equivalent with $\Lambda'$, and get
\begin{prop} \label{conditiiledecoupling2} If the pair
$(\Lambda,E)$ is $\mathcal{F}$-coupling of the second kind, it is
a Jacobi structure iff $\Lambda''$ is a leaf-tangent Poisson
structure, {\rm(\ref{LieinJac2})} hold and
\begin{equation} \label{VorobJac2} d'\sigma=
(\flat_\sigma
E)\wedge\sigma,\,L_X\Lambda''=\sigma(X,E)\Lambda'',\;
-p_F[X,Y]=\sharp_{\Lambda''}\{d[\sigma(X,Y)]\},
\end{equation} for all the vector fields
$X,Y\in\mathcal{V}^{1,0}_{pr}(M).$ \end{prop}

On the other hand we have \begin{prop}
\label{conditiiledecoupling3} The $\mathcal{F}$-coupling of the
third kind pair $(\Lambda,E)$ is a Jacobi structure iff $M$ is a
locally product manifold between the foliation $\mathcal{F}$ and a
foliation $\mathcal{H}$ and, with respect to this locally product
structure, $\Lambda''=0$ and $(\Lambda',E)$ is a Jacobi structure
on $M$ such that its restrictions to the leaves of $\mathcal{H}$
are contact structures of these leaves.
\end{prop} \noindent {\bf Proof.} The third kind coupling
conditions imply the existence of a $(1,0)$-form $\xi\in
ker\,\sharp_{\Lambda'}$ such that $\xi(E)=1$. Then, if we take
$\alpha=\xi$ in the third condition (\ref{farabarJac2}), we get
$\Lambda''=0$, and we remain with the first condition
(\ref{farabarJac2}) and with
$$\lambda([\sharp_{\Lambda'}\alpha,\sharp_{\Lambda'}\beta])=0,\;
L_E\Lambda'=0.$$ This implies the integrability of the normal
distribution $H$, which provides the foliation $\mathcal{H}$, and
shows that $(\Lambda',E)$ is the required Jacobi structure. Q.e.d.
\begin{rem}\label{obsvecinfoicontact}
{\rm Proposition \ref{conditiiledecoupling3} may be applied to a
tubular neighborhood of a contact characteristic leaf of an
arbitrary Jacobi manifold $(M,\Lambda,E)$ (see Example
\ref{exfoiinJacobi}) and it tells us that all the neighboring
characteristic leaves also are contact manifolds. This fact may be
seen straightforwardly as follows. The characteristic leaf
$L_{x_0}$ through the point $x_0\in M$ is a contact leaf iff there
exists a $1$-form $\theta\in ann(im\,\sharp_\Lambda)$ such that
$\theta_{x_0} (E_{x_0})\neq0$, and, if this condition holds at
$x_0$, it holds on an open neighborhood of $x_0$. The situation is
different for a locally conformal symplectic leaf $L_{x_0}$. For
instance, take $M=\mathbb{R}^4=\{(q,p,u,t)\}$,
$$\Lambda=u\frac{\partial}{\partial q}\wedge \frac{\partial}{\partial p}
+ (t\frac{\partial}{\partial t})\wedge(p\frac{\partial}{\partial
p}),\;\;E=t\frac{\partial}{\partial t},$$ and
$x_0(q=1,p=0,u=1,t=0)$. Then, $L_{x_0}$ is $2$-dimensional but,
any neighboring leaf $L_{x_1}$ where $x_1(q=1,p=0,u=1,t\neq0)$ is
$3$-dimensional, hence, a contact leaf.}\end{rem}

It is interesting to discuss the coupling conditions for the {\it
transitive} Jacobi structures. i.e., locally conformal symplectic
structures and contact structures \cite{{DLM},{V3}}.

Recall that a locally conformal symplectic structure of a
differentiable manifold $M$ is a non degenerate $2$-form $\omega$,
which satisfies the condition
\begin{equation} \label{condlcs} d\omega=\epsilon\wedge\omega,
\end{equation} where $\epsilon$ is a closed $1$-form (the {\it Lee
form}). Then, the bivector field $\Lambda$ defined by
$\sharp_\Lambda=-\flat_\omega^{-1}$ and the vector field
$E=\sharp_\Lambda\epsilon$ define a Jacobi structure on $M$. In
the general case of a manifold $M$ endowed with a pair
$(\omega,\epsilon)$ where $\omega$ is a non degenerate  $2$-form
and $\epsilon$ is a $1$-form, a pair $(\Lambda,E)$ can be defined
similarly, and it may be coupling of either the first or the
second kind for a foliation $\mathcal{F}$. This happens iff the
pullback of $\omega$ to $\mathcal{F}$ is non degenerate, and then
$H=\sharp_\Lambda(ann\,F)$ is the $\omega$ orthogonal distribution
of $F$, and $E$ is always tangent or nowhere tangent to
$\mathcal{F}$, respectively. In these cases, we also attribute the
coupling property to the pair $(\omega,\epsilon)$ itself.
Accordingly, if we want the coupling property in the locally
conformal symplectic case, $\omega$ must be of the form
(\ref{descsymplcoup}) with respect to the bigrading defined by the
symplectic-orthogonal decomposition $TM=H\oplus F$ and satisfy
(\ref{condlcs}). The result is
\begin{prop}
\label{lcscoupling} Let $$\omega=\sigma_{2,0}+\theta_{0,2}$$ be a
non degenerate $\mathcal{F}$-coupling $2$-form on
$(M^{2n},\mathcal{F}^{2s})$. If the coupling is of the first kind
$(E\in\Gamma F)$, $\omega$ is locally conformal symplectic iff
\begin{equation}\label{eqlcscoup1}
d'\sigma=0,\;d'\theta=0,\;d''\sigma+\partial\theta=\epsilon
\wedge\sigma,\;d''\theta=\epsilon\wedge\theta. \end{equation}
If the coupling is of the second kind, $\omega$ is locally
conformal symplectic iff
\begin{equation}\label{eqlcscoup2}
d'\sigma=\epsilon\wedge\sigma,\;d''\sigma+\partial\theta=0,\;
d'\theta=\epsilon\wedge\theta,\; d''\theta=0.
\end{equation}
\end{prop} \noindent{\bf Proof.} The Lee form $\epsilon$ has
bidegree $(0,1)$ in the first case, and $(1,0)$ in the second. The
results follow by equating the corresponding homogeneous
components of the two sides of (\ref{condlcs}). Q.e.d.

In the first kind coupling case, the leaves of $\mathcal{F}$ get
the induced locally conformal symplectic structure
$(\theta,\epsilon)$ and in the second kind coupling case the
leaves get the induced symplectic structure $\theta$. If the forms
$(\epsilon_{0,1},\theta_{0,2})$, which define leafwise locally
conformal symplectic structures of the leaves of $\mathcal{F}$,
and the normal bundle of $\mathcal{F}$ are given, and we look for
the form $\sigma$ that extends the given data to a locally
conformal symplectic form $\omega$, we may encounter obstructions.
In particular, $(M,\epsilon)$ $(d\epsilon=0)$ has a {\it twisted}
de Rham cohomology $H^*_\epsilon(M)$, which is that of the cochain
complex $(\Omega^*(M),d_\epsilon=d-\epsilon\wedge)$. Using the
degree decomposition of $d^2=0$, we see that the form
$\partial\theta$ is $d_\epsilon$-closed, and the cohomology class
$[\partial\theta]\in H^3_\epsilon(M)$ is an obstruction to the
existence of $\sigma$. The same holds for a given pair
$(\epsilon_{1,0},\theta_{0,2})$ which satisfies the last two
conditions (\ref{eqlcscoup2}) and which we would like to extend to
a locally conformal form $\omega$ on $M$.

A contact form $\phi$ on a differentiable manifold $M^{2n+1}$ is a
$1$-form such that $\phi\wedge (d\phi)^n\neq0$ at each point of
$M$. The contact form $\phi$ defines the field of {\it contact
hyperplanes} $\pi$ of equation $\phi=0$ and the {\it Reeb vector
field} $E$ characterized by
\begin{equation}\label{Reeb} \phi(E)=1,\;
i(E)d\phi=0,\end{equation} and we may write
\begin{equation}\label{descReeb} TM=\pi\oplus span\{E\},\;
T^*M=\pi^*\oplus span\{\phi\}.\end{equation} The mapping
$\flat_{d\phi}$ is an isomorphism of $\pi$ onto $\pi^*$, and we
may define a bivector field $\Lambda$ by asking $\sharp_\Lambda$
to be zero on $span\{\phi\}$ and $-\flat_{d\phi}^{-1}$ on $\pi^*$.
It turns out that the pair $(\Lambda,E)$ is a Jacobi structure
\cite{{DLM},{V3}}.

Let $(M,\phi)$ be a contact manifold. We will look at two kinds of
foliations $\mathcal{F}$ of $M$: {\it foliations of the contact
type}, characterized by the condition $E\in\Gamma F$
$(F=T\mathcal{F})$, and {\it foliations of the symplectic type},
characterized by $F\subseteq\pi$.

If $\mathcal{F}$ is of the contact type, only coupling of the
first kind can occur. More exactly, if $E\in\Gamma F$, $F$ and
$\pi$ are transversal, $F\cap\pi$ is a vector subbundle of $TM$,
$ann\,F\subseteq\pi^*$ and
$\flat_{d\phi}^{-1}(ann\,F)=(F\cap\pi)^{\perp_{d\phi}}$. The
coupling condition holds iff $d\phi$ is non degenerate on
$(F\cap\pi)^{\perp_{d\phi}}$, equivalently, $d\phi$ is non
degenerate on $(F\cap\pi)$. This means that the Jacobi structure
defined by the contact form $\phi$ is coupling for $\mathcal{F}$
iff the pullbacks of $\phi$ to the leaves of $\mathcal{F}$ are
contact forms of the leaves. The corresponding normal bundle $H$
of $\mathcal{F}$ is $(F\cap\pi)^{\perp_{d\phi}}$. The problem of
extending a smooth family of contact structures of the leaves of
$\mathcal{F}$ to a coupling contact structure of $M$ (in the sense
of Jacobi structures) was recently studied for fiber bundles by
Lerman \cite{Lr}.

If $\mathcal{F}$ is of the symplectic type, $ann(F \oplus E)=
(ann\,F)\cap\pi^*$ and $\flat_{d\phi}^{-1}[(ann\,F)\cap\pi^*]=
F^{\perp_{d\phi}}\subseteq\pi$, which does not contain $E$.
Accordingly, only couplings of the third kind might occur, and the
Jacobi structure defined by the contact form $\phi$ would be
coupling for $\mathcal{F}$ iff the pullbacks of $d\phi$ to the
leaves of $\mathcal{F}$ would be symplectic forms of the leaves.
This is impossible since $d\phi=0$ along any integral submanifold
of the contact distribution $\pi$.

We end this paper by indicating some possible applications of the
Vorobiev-Poisson structures (Section 4) to Jacobi manifolds.

First we recall that, if $(M,\Lambda,E)$ is a Jacobi manifold, the
jet bundle $J^1(M,\mathbb{R})\approx T^*M\oplus
\mathbb{R}$ has a natural structure of a Lie algebroid, induced
by the restriction to $M\approx M\times\{0\}$ of the cotangent Lie
algebroid of the Poisson manifold $(M\times\mathbb{R},P)$ with $P$
given by (\ref{homogPsptJacobi}) (e.g., see \cite{V3}). The
corresponding Lie bracket is
\begin{equation}
\label{crosetLieptJacobi} \{(\alpha,f),(\beta,g)\}=
(\{\alpha,\beta\}_\Lambda+fL_E\beta-gL_E\alpha \end{equation}
$$-\alpha(E)\beta+\beta(E)\alpha,
\{f,g\}-\Lambda(df-\alpha,dg-\beta)),
$$ where $\alpha,\beta\in\Omega^1(M),
f,g\in C^\infty(M)$, the bracket of functions is (\ref{Jacobibr})
and the $\Lambda$-bracket of $1$-forms is (\ref{bracket1}). the
The corresponding anchor map is
\begin{equation} \label{ancorapeS}
\rho(\alpha,f)=\sharp_\Lambda\alpha+fE. \end{equation}

Since the characteristic leaves of the Jacobi structure are
integral submanifolds of the distribution
$span\{im\,\sharp_\Lambda,\,E\}$, the bracket
(\ref{crosetLieptJacobi}) actually computes along these leaves
and, for any characteristic leaf $S$ of $(\Lambda,E)$,
$J^1(M,\mathbb{R})|_S$ is a transitive Lie algebroid with the
basis $S$.

Furthermore, we have
\begin{prop}\label{nucleuptJacobi} For any characteristic leaf
$S$, the kernel $K=ker\,\rho$ of the anchor map
{\rm(\ref{ancorapeS})} restricted to $S$ is isomorphic with the
annihilator $ann_{M\times\mathbb{R}}(TS)$ of $S$ seen as the
submanifold $S\times\{0\}$ of $M\times\mathbb{R}$. \end{prop}
\noindent {\bf Proof.} If the dimension of $S$ is even $S$ is a
locally conformal symplectic manifold,
$TS=im\,\sharp_{\Lambda|_S}$, and there exist $1$-forms
$\epsilon\in T^*M|_S$ such that $E|_S=\sharp_\Lambda(\epsilon)$.
As in formula (\ref{homogPsptJacobi}), we denote by $t$ the
coordinate on $\mathbb{R}$. The mapping $(\alpha,f)\mapsto
(\alpha+f\epsilon) +fdt$ yields an isomorphism $K\rightarrow
ann_{M\times\mathbb{R}}(TS)$ as required, which is not unique
since it depends on the choice of $\epsilon$.

If the dimension of $S$ is odd $S$ is a contact manifold with the
Reeb vector field $E|_S$, which neither vanishes nor belongs to
$im\,\sharp_\Lambda$ at each point of $S$. Accordingly,
{\rm(\ref{ancorapeS})} yields $K=ker\,\sharp_\Lambda$ along $S$.
On the other hand, $ann_M(TS)=(ker\,\sharp_\Lambda)\cap ann(E)$.
Hence, if we choose the $1$-form
$\xi\in\Gamma(ker\,\sharp_\Lambda)$ such that $\xi(E)=1$ along
$S$, we get an isomorphism $K\approx ann_M(TS)\oplus span\{\xi\}$,
which allows us to represent the elements of $K$ under the form
$\beta+h\xi$ where $\beta\in ann_M(TS),h\in C^\infty(M)$.
Furthermore, the mapping $(\beta+h\xi,0)\mapsto \beta+hdt$ yields
an isomorphism $K\rightarrow ann_{M\times\mathbb{R}}(TS)$ as
required. Again this isomorphism is not unique since it depends on
the choice of $\xi$. Q.e.d.

Now, assume that $S$ is an embedded characteristic leaf of the
Jacobi manifold $(M,\Lambda,E)$ and that the induced structure of
$S$ is globally conformal symplectic, i.e., there exists on $S$ a
symplectic form $\omega$ and a function $f\in C^\infty(M)$ such
that the Jacobi structure of $S$ is defined by the $2$-form
$e^f\omega$. The symplectic form $\omega$ allows us to define a
Vorobiev-Poisson structure $\Pi$ on a neighborhood $U$ of $S$ seen
as the zero section of the total space of the Lie coalgebras
bundle $K^*$. By Proposition \ref{nucleuptJacobi}, $U$ may also be
seen as a tubular neighborhood of $S\times\{0\}$ defined by a
normal bundle of this submanifold in $M\times\{\mathbb{R}\}$.
Using formula (\ref{conformJ}) with $a=e^{f},\Lambda=\Pi,E=0$ we
get a Jacobi structure on $U$ for which $S$ is a characteristic
leaf with the same induced Jacobi structure as the one induced by
the original $(\Lambda,E)$. Notice that we cannot just use
Vorobiev's construction along a locally conformal symplectic leaf
since the form (\ref{sigmapealgLie}) associated to a locally
conformal symplectic form $\omega$ does not satisfy conditions
(\ref{VorobJac2}).

On the other hand, if $S$ is a contact, embedded, characteristic
leaf of the Jacobi manifold $(M,\Lambda,E)$,
$S\times\{\mathbb{R}\}$ is a symplectic leaf of the Poisson
structure $P$ of $M\times\{\mathbb{R}\}$ and its symplectic
structure is the so-called {\it symplectification} of the contact
structure of $S$. Using this symplectic structure, we get a
Vorobiev-Poisson structure on a neighborhood of the zero section
of a normal bundle $N(S\times\{\mathbb{R}\})$ in
$M\times\{\mathbb{R}\}$. This neighborhood may be seen as having
the form $V\times\{\mathbb{R}\}$ where $V$ is a neighborhood of
$S$ seen as the zero section of a normal bundle $NS$ of $S$ in
$M$. Therefore, $S$ has a neighborhood $V\subseteq NS$ which is a
hypersurface of a Vorobiev-Poisson manifold.
\begin{rem}\label{obsfinala} {\rm If $L_{x_0}$ is a locally conformal
symplectic leaf of $(M,\Lambda,E)$ with the $2$-form $\omega$ and
the (closed) Lee form $\epsilon$, the symplectic leaf of
$(M\times\mathbb{R},P)$ through $(x_0,0)$ is the hypersurface of
$(L_{x_0}\times\mathbb{R},P)$, which integrates the equation
$\epsilon-dt=0$, and passes through $x_0$.}\end{rem}
 %\end{center}
\hspace*{7.5cm}{\small \begin{tabular}{l} Department of
Mathematics\\ University of Haifa, Israel\\ E-mail:
vaisman@math.haifa.ac.il \end{tabular}}

\begin{thebibliography}{xx}
\bibitem{GD} I. M. Gelfand and I. Ya. Dorfman, The Schouten
bracket and Hamiltonian operators. Funkt. Anal. Prilozhen. 14 (3)
(1980), 71-74.
\bibitem{DLM} P. Dazord, A. Lichnerowicz and Ch.-M. Marle,
Structures locales des vari\'et\'es de Jacobi, J. Math. pures et
appl., 70 (1991), 101-152.
\bibitem{GLS} V. Guillemin, E. Lerman and S. Sternberg, Symplectic
fibrations and multiplicity diagrams. Cambridge Univ. Press,
Cambridge, 1996.
\bibitem{MIgl} D. Iglesias and J. C. Marrero, Some linear Jacobi
structures on vector bundles, C. R. Acad. Sci. Paris, S\'er. I
Math., 331  (2000), 125-130.
\bibitem{Lr} E. Lerman, Contact fiber bundles,
arXiv:math.DG/0301137.
\bibitem{Lch} A. Lichnerowicz, Global theory of connections and
holonomy groups. Noordhoof Intern. Publ., Leyden, 1976.
\bibitem{Mac} K. C. H. Mackenzie, Lie Groupoids and Lie Algebroids
in Differential Geometry, LMS Lecture Notes Ser., Vol. 124,
Cambridge Univ. Press, Cambridge, 1987.
\bibitem{MX} K. C. H. Mackenzie and P. Xu, Lie bialgebroids and
Poisson groupoids, Duke Math. J., 18 (1994), 415-452.
\bibitem{Mol} P. Molino, Riemannian foliations. Progress in Math.
73, Birkh\"auser, Boston, 1988.
\bibitem{MMR} R. Montgomery, J. E. Marsden and T. Ra\c tiu, Gauged
Lie-Poisson structures, Cont. Math. AMS, Vol. 28 (Boulder
Proceedings on Fluids and Plasmas), 1984, 101-114.
\bibitem{St} S. Sternberg, On minimal coupling and the symplectic
mechanics of a classical particle in the presence of a Yang-Mills
field, Proc. Nat. Acad. Sci. USA, 74 (1977), 5253-5254.
\bibitem{V4} I. Vaisman, Cohomology and differential forms. M.
Dekker, Inc., New York, 1973.
\bibitem{V1} I. Vaisman, Lectures on
the Geometry of Poisson Manifolds. Progress in Math. 118,
Birkh\"auser, Basel, 1994.
\bibitem{V3} I. Vaisman, A lecture on Jacobi manifolds. Selected
topics in Geom. and Math. Physics, 1 (2002), 81-100.
\bibitem{V5} I. Vaisman, Hamiltonian structures on foliations, J.
Math. Physics, 43 (2002), 4966-4977.
\bibitem{Vor} Y. Vorobjev, Coupling tensors and Poisson geometry
near a single symplectic leaf. In: Lie algebroids and related
topics in differential geometry (J. Kubarski, P. Urba\'nski and R.
Wolak, eds.), Banach Center Publ., Vol. 54, Warszawa, 2001, p.
249-274.
\end{thebibliography}
\end{document}